\documentclass[11pt]{article}

\usepackage{tikz}
\usepackage{subfigure}
\usepackage[english]{babel}
\usepackage{graphicx}
\usepackage{thm-restate}

\usepackage[center]{caption2}
\usepackage{amsfonts,amssymb,amsmath,latexsym,amsthm}
\usepackage{multirow}

\def\diam{\mbox{diam}}
\def\dist{\mbox{dist}}

\newtheorem{thm}{Theorem}[section]
\newtheorem{defi}[thm]{Definition}
\newtheorem{cor}[thm]{Corollary}
\newtheorem{lem}[thm]{Lemma}
\newtheorem{prop}{Proposition}

\newtheorem{conj}[thm]{Conjecture}
\newtheorem{prob}[thm]{Problem}

\newtheorem{claim}{Claim}

\newtheorem{case}{Case}%[section]

\begin{document}

\title{Dirac-type condition for Hamilton-generated graphs
%\thanks{The work was supported by NNSF of China (No. 12071453) and  Anhui Initiative in Quantum Information Technologies (AHY150200) and National Key Research and Development Project (SQ2020YFA070080).}
}
\author{Xinmin Hou$^{a,b}$, \quad Zhi Yin$^{a}$\\
%\small $^{a}$ Key Laboratory of Wu Wen-Tsun Mathematics\\
\small $^{a}$ School of Mathematical Sciences\\
\small University of Science and Technology of China, Hefei 230026, Anhui, China\\
\small$^b$ Hefei National Laboratory\\
\small University of Science and Technology of China,  Hefei 230088, Anhui, China\\
\small Email: xmhou@ustc.edu.cn (X. Hou), yinzhi@mail.ustc.edu.cn (Z. Yin)
}

\date{}

\maketitle

\begin{abstract}
%We study Hamilton-generated graphs under Dirac-type conditions. 
The cycle space $\mathcal{C}(G)$ of a graph $G$ is defined as the linear space spanned by all cycles in $G$. For an integer $k\ge 3$, let $\mathcal{C}_k (G)$ denote the subspace of $\mathcal{C}(G)$ generated by the cycles of length exactly   $k$.
A graph $G$ on $n$ vertices is called Hamilton-generated if $\mathcal{C}_n (G) = \mathcal{C}(G)$, meaning every cycle in $G$ is a symmetric difference of some Hamilton cycles of $G$.
%A necessary condition for this property is that $n$ must be odd. 
Heinig (European J. Combin., 2014) showed that for any $\sigma >0$ and  sufficiently large	odd $n$, every $n$-vertex graph with minimum degree $(1+ \sigma)n/2$ is Hamilton-generated. He further posed the question that whether the minimum degree requirement  could be lowered to the Dirac threshold $n/2$. Recent progress by Christoph, Nenadov, and Petrova~(arXiv:2402.01447) reduced the minimum degree condition to $n/2 + C$ for some large constant $C$. 
In this paper, we resolve Heinig's problem completely by proving that for sufficiently large odd $n$, every Hamilton-connected graph $G$ on $n$ vertices with minimum degree at least $(n-1)/2$ is Hamilton-generated. 
Moreover, this result is tight for the minimum degree and the Hamilton-connected condition.
The proof relies on the parity-switcher technique introduced by Christoph, et al in their recent work, as well as a classification lemma that strengthens a previous result by Krivelevich, Lee, and Sudakov~(Trans. Amer. Math. Soc., 2014).

\end{abstract}

\section{Introduction}

For a graph $G=(V,E)$, a subset  of $E$  corresponds naturally to a vector (i.e., its indicator function $f : E\rightarrow \mathbb{F}_2$) in $\mathbb{F}_2^{|E|}$.
For convenience, we do not distinguish a subset $F$ in $E$ and its corresponding vector in $\mathbb{F}_2^{|E|}$.
 The {\em edge space}  $\mathcal{E}(G)$  of $G$ is the vector space over the field $\mathbb{F}_2$ formed by the subsets of $E$.
The {\em cycle space} $\mathcal{C}(G)$ is the subspace of $\mathcal{E}(G)$ spanned by all cycles in $G$. 
Denote $\mathcal{C}_k (G)$ as the subspace of $\mathcal{C}(G)$ spanned by the cycles of length exactly  $k$ in $G$.
A graph $G$ on $n$ vertices is said to be {\em Hamilton-generated} if $\mathcal{C}_n (G) = \mathcal{C}(G)$, and otherwise, it is called  {\em non-Hamilton-generated}.
A {\em Hamilton path} is a path that visits each vertex in a graph exactly once. Additionally, a graph is {\em Hamilton-connected} if every pair of vertices is connected by a Hamilton path. 

A natural question is to determine the condition under which $\sum_{k\in I} \mathcal{C}_k(G) =\mathcal{C}(G)$ for some finite positive integer set $I$. 
According to \cite{hartman1983}, a conjecture originating with Bondy states that:
%, is a well-studied problem in graph theory \cite{Baron,Kahn,baron2019,bondy1981,demarco2013,hartman1983}, and in 1979, J. A. Bondy has following conjecture, which is still open now.
\begin{conj} [Bondy, 1979]\label{CONJ:Bondy}
   For every $d \in \mathbb{Z}$, in every $3$-connected graph $G$ with
   $|V(G)| \ge 2d$ and $\delta(G) \ge d$, $$\sum_{k\ge 2d-1} \mathcal{C}_k(G) =\mathcal{C}(G).$$
   %any cycle in $G$ can be represented by the set of all cycles of length at least $2d-1$ in $\mathbb{F}_2$.
\end{conj}
This conjecture is still open. In the case that $G$ is 2-connected, Hartman~\cite{hartman1983} proved (a variant of Conjecture~\ref{CONJ:Bondy}) that  
$\sum_{k\ge d+1} \mathcal{C}_k(G) =\mathcal{C}(G)$ for each 2-connected graph $G$ with minimum degree $d$ and $G\not=K_{d+1}$ if $d$ is odd.
Locke (\cite{locke85}, Theorem 4) proved  that Bondy’s conjecture is
true under the additional assumption of ‘$G$ non-Hamiltonian or $|V(G)| > 4d-5$’. For Hamiltonian graphs, Barovich, Locke~\cite{BC00}
showed that Bondy’s conjecture is true for a large class of Hamiltonian graphs. Determining conditions under which $\mathcal{C}_k(G) = \mathcal{C}(G)$, for some $k$, is a well-studied problem in graph theory \cite{baron2019,bondy1981,demarco2013,hartman1983,locke85a,locke85}.

A classical theorem of Dirac from 1952 asserts that every graph on $n$ vertices with minimum degree at least $\lceil n/2\rceil$ is Hamiltonian.
In this paper, we primarily focus on  Bondy’s conjecture under Dirac-type conditions when the size of $G$ is large.
We call a graph $G$ on $n$ vertices {\em Hamilton-generated} if  $\mathcal{C}_n (G) = \mathcal{C}(G)$.
Note that if $G$ contains an odd cycle, a necessary condition for $\mathcal{C}_n (G) = \mathcal{C}(G)$ is that $G$ has an odd number of vertices. 
Heinig~\cite{Heinig} established the following result.
\begin{thm}[Heinig, \cite{Heinig}]\label{THM: Heinig}
For any $\sigma >0$ and  sufficiently large	odd $n$, if $G$ is a graph on $n$ vertices with minimum degree $(1+ \sigma)n/2$, then $\mathcal{C}_n (G) = \mathcal{C}(G)$.
\end{thm} 
 Heinig further asked the following question in~\cite{Heinig}.
\begin{prob} [\cite{Heinig}] \label{Prob}
     Suppose $n$ is odd and large enough, and let $G$ be a graph with $n$ vertices and $\delta(G)> n/2$, is it true that  $\mathcal{C}_{n}(G)=\mathcal{C}(G)$?
\end{prob}
A recent progress  by Christoph, Nenadov, and Petrova~\cite{Nenadov} lowered the minimum degree condition to $n/2 + C$ for some large constant $C$.
%by providing  an almost optimal lower bound for the minimum degree,  up to an additive constant term.
\begin{thm}[\cite{Nenadov}]\label{Nenadov}
    There exists a constant $C>1$ such that the following holds for sufficiently large odd $n$. Let $G$ be a graph with $n$ vertices and $\delta(G) \ge n/2 + C$. Then $\mathcal{C}_{n}(G)=\mathcal{C}(G)$.
\end{thm}

In this paper, we completely resolve Problem~\ref{Prob} affirmatively.
In fact, we present a tight lower bound for the minimum degree required to insure a graph is Hamilton-generated,  as demonstrated below. 

\begin{restatable}{thm}{Odd} \label{Odd}
    There exists \( n_0 \in \mathbb{N}^+ \) such that for any odd \( n > n_0 \), if \( G \) is a graph with \( n \) vertices that has a minimum degree \( \delta(G) \ge (n-1)/2 \)  and is Hamilton-connected, then \( \mathcal{C}_{n}(G) = \mathcal{C}(G) \).
\end{restatable}

Note that any graph $G$ with $\delta(G)> |V(G)|/2$ is necessarily Hamilton-connected. Consequently,  Problem~\ref{Prob} is fully resolved as a corollary of Theorem~\ref{Odd} as demonstrated below.
\begin{cor}
If $G$ is a graph with odd $n$ vertices and $\delta(G)> n/2$, then $\mathcal{C}_{n}(G)=\mathcal{C}(G)$ for sufficiently large $n$.  
\end{cor}

The following constructions show that Theorem \ref{Odd} is tight for the minimum degree and the Hamilton-connected condition.

\noindent{\bf Construction A:} Let $n=4k+1$ for some $k \ge 2$, and let $X$ and $Y$ be two disjoint sets with  $|X|=(n+1)/2=2k+1$ and $|Y|=(n-1)/2=2k$. Suppose $a_1,a_2,a_3 \in X$, $b_1,b_2,b_3 \in Y$.  Define graphs $G_1, G_2$ and $G_3$  on the vertex set $V=X \cup Y$ with the following edge sets: 
\begin{align*}
 E(G_1)&=\{xy: x \in X, y \in Y\} \cup \{a_1a_2,b_1b_2\},  \\
 E(G_2)&=\{xx':x,x' \in X \cup \{b_1\}\} \cup \{yy':y,y' \in Y \cup \{a_1\}\} \setminus \{a_1b_1\},\\
 E(G_3)&=\{xx':x,x' \in X\} \cup \{yy':y,y' \in Y\} \cup \{a_1b_1,a_2b_2,a_3b_3\}.
\end{align*}
In other words, $E(G_1)=E(K[X,Y])\cup\{a_1a_2, b_1b_2\},\, E(G_2)=(E(K[X\cup\{b_1\}])\cup E(K[Y\cup\{a_1\}])) \setminus \{a_1b_1\}$, and $E(G_3)=E(K[X])\cup E(K[Y])\cup \{a_1b_1, a_2b_2, a_3b_3\}$, where $K[U, V]$ denotes the complete bipartite graph with bipartite sets $U$ and $V$, and $K[V]$ denotes the complete graph with vertex set $V$.  
\begin{figure}[htbp]
    \centering
    \includegraphics[width=1\textwidth]{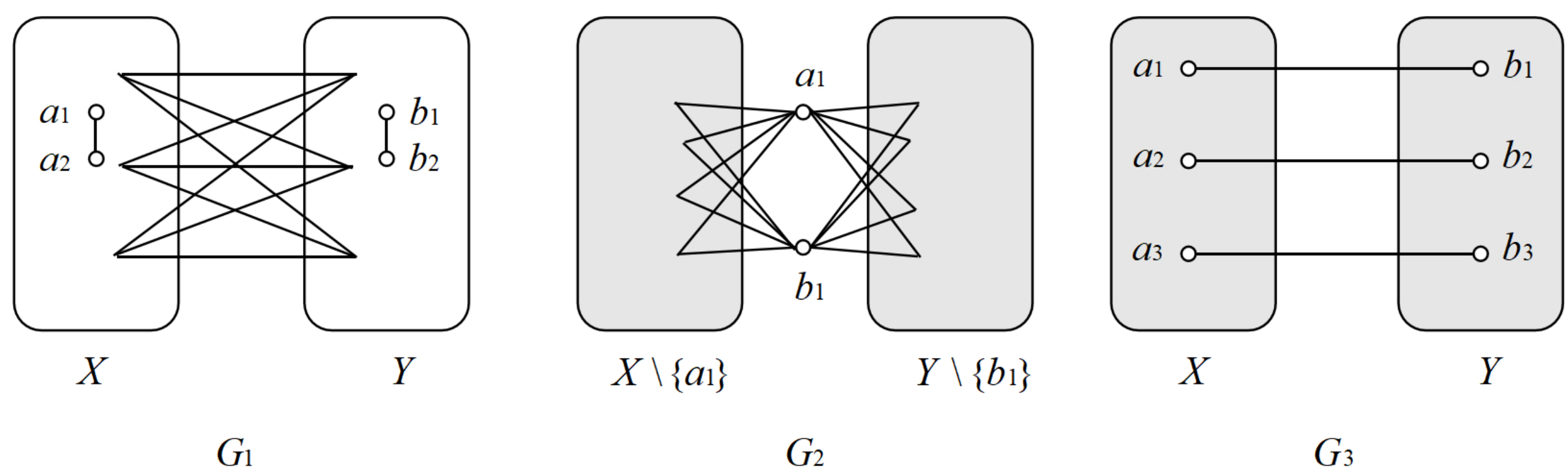}  % 图片路径
    \caption{$G_1,G_2,G_3$ in Construction A}
    \label{fig:example}
\end{figure}

From Construction A, we can check that $\delta(G_1)=\delta(G_2)=\frac{n-1}2$. Both $G_1$ and $G_2$  contain Hamilton  cycles, but  they are  not Hamilton-connected. Additionally, $\delta(G_3) =\frac{n-3}2$ and it is Hamilton-connected. However,   $G_1,G_2$ and $G_3$ are not Hamilton-generated (Proposition~\ref{example}), implying that  the minimum degree and Hamilton-connectedness are both optimal in a sense.

\begin{prop}\label{example}
 Graphs   $G_1,G_2$ and $G_3$ are not Hamilton-generated.
\end{prop}
\begin{proof}
    For $G_1$, it is easy to see that every Hamilton cycle in $G_1$ does not contain the edge $b_1b_2$.
   Hence, any cycles containing $b_1b_2$ cannot be generated by Hamilton cycles (in fact, there are numerous such cycles containing $b_1b_2$ in $G_1$). 
   
   For $G_2$, let $X'=X\setminus \{a_1\}$ and $Y'=Y\setminus \{b_1\}$. 
   Then every Hamilton cycle contains
    precisely one edge in $E(a_1,X')$, in $E(a_1,Y')$, in $E(b_1,X')$, and in $E(b_1,Y')$, where $E(X,Y)=\{xy: x \in X, y \in Y\}$ for subsets $X,Y$ in $V(G)$. Hence
    the symmetric difference of an odd number of Hamilton cycles must contain at least one edge from each of $E(a_1,X')$ and $E(a_1,Y')$. 
   Note that any odd cycle should be represented by an odd number of Hamilton cycles. Therefore, any odd cycle not containing $a_1$ can not be represented by Hamilton cycles.
    
  For $G_3$, since there are exactly three edges between $X$ and $Y$, it is easy to check that any Hamilton cycle of $G$ contains a path $P$ such that $P$ is a Hamilton path of $G[X]$. Note that $|X|=2k+1$. Hence $P$ has an even length. This means that any Hamilton cycle contains an even number of edges in $G[X]$.  Therefore, any element in $\mathcal{C}_n(G)$ also contains an even number of edges in $G[X]$.  This implies that $\mathcal{C}_n(G)\not=\mathcal{C}(G)$ as any odd cycle in $G[X]$ would not belong to $\mathcal{C}_n(G)$.
\end{proof}

\section{Tools and Lemmas}

We first introduce the additional notation used throughout the paper.
For a graph $G=(V,E)$ and subsets $X,Y \subseteq V$, let $E_G(X)=E(G[X])$, $E_G(X,Y)=\{xy: x \in X, y \in Y\}$, and $G[X,Y]$ be the bipartite graph with vertex parts $X \cup Y$ and edge set $E_G(X,Y)$. Denote $e_G(X)=|E_G(X)|$ and $e_G(X,Y)=|E_G(X,Y)|$ (the subscript will be omitted  if the graph $G$ is clear from the context). 
For convenience, write $E_G(x,Y)=E_G(\{x\},Y)$ and $e_G(x,Y)=e_G(\{x\},Y)$.
Define  $\delta_X(Y)=\min\{|N_G(y) \cap X|:y \in Y\}$ and $\Delta_X(Y)=\max\{|N_G(y) \cap X|:y \in Y\}$. 
We say $\mathcal{P}$ is \textit{an union of disjoint paths} in $G$ if $\mathcal{P}$ is a subgraph of $G$ such that $\mathcal{P}$ is a forest and $1 \le \Delta(\mathcal{P}) \le 2$, denoted as \textit{UDP}.
If $\mathcal{P}$ is a UDP in $G$, we denote the set of endpoints of paths in $\mathcal{P}$ as End$(\mathcal{P})$, and let In$(\mathcal{P})=V(\mathcal{P})\setminus \text{End}(\mathcal{P})$. Additionally, denote the number of paths in $\mathcal{P}$ by $|\mathcal{P}|$. 
We observe that $$|\text{End}(\mathcal{P})|+|\text{In}(\mathcal{P})|=|\text{In}(\mathcal{P})|+2|\mathcal{P}|=|V(\mathcal{P})|=|E(\mathcal{P})|+|\mathcal{P}|.$$

For subsets  $W\subseteq V$ and $E_0 \subseteq \tbinom{V}{2}$, let $G- W=G[V\setminus W]$ and $G+E_0$ be the graph with  vertex set $V$ and edge set $E \cup E_0$. If $R$ is a subgraph of $G$, we let $G\setminus R$ be the spanning subgraph of $G$ induced by the edge set $E(G)\setminus E(R)$.
For two sets $A$ and $B$, write $A\setminus B$ for $A\setminus (A\cap B)$.

 The following lemma is proposed in \cite{Nenadov} and a similar statement is also used as a starting point in \cite{baron2019,demarco2013}.

\begin{lem} [\cite{Nenadov}] \label{R-lem} 
   Let $G$ be a Hamiltonian graph with an odd number $n$ of vertices such that $\mathcal{C}_{n}(G) \ne \mathcal{C}(G)$.  Then there
   exists a subgraph $R \subseteq G$ for which the following holds$:$
   \begin{itemize}
       \item[(a)] $R \ne G$,
       \item[(b)] every Hamilton cycle in G contains an even number of edges from R, and
       \item[(c)] for every partition $V (G) = A \cup B$, we have $e_{R}(A,B) \ge \frac 12e_{G}(A,B)$ and $R \ne G[A,B]$.
   \end{itemize}
\end{lem}
%{\color{blue}这个lemma的序号看起来有点奇怪。}

%Next, we define the Parity-switcher, which is an important subgraph of the non-Hamilton-generated graph.

\begin{defi}[Parity-switcher]
    Given a graph $G$ and a subgraph $R \subseteq G$, a subgraph $W \subseteq G$ is called
    an \textbf{$R$-parity-switcher} if it consists of an even cycle $C = (v_1, v_2, ... , v_{2k}, v_1)$ with an odd number of edges in
    $R$, and vertex-disjoint paths $P_i$ for $1\le i\le k+1$ with $P_i$ connects $v_{i}$ and $v_{2k-i+2}$  for $2 \le i \le k$ and $P_j$ connects $v_j$ and $u_j$ for $j\in\{1,k+1\}$, where $u_1$ and $u_{k+1}$ are vertices in $G$ and referred to as \textbf{touch-vertices}.
    Denote such an $R$-parity-switcher as $W=W(R, C_{(v_1, v_2, ... , v_{2k})}, P_1,P_2,\ldots, P_k, P_{k+1}, u_1, u_{k+1})$.

%Additionally, for an $R$-parity switcher \(W\), a subgraph $W'=(W,Q_1,Q_{k+1})$ extended from $W$ is called an \textbf{extended \( R \)-parity-switcher} if it consists of \( W \) and two vertex-disjoint paths \( Q_1 \) and \( Q_{k+1} \) such that
%\( V(W) \cap V(Q_i) = \{v_i\} \) for $i\in\{1, k+1\}$.
%In this case, we refer to the two endpoints of \( Q_1 \) and \( Q_{k+1} \), other than \( v_1 \) and \( v_{k+1} \), as the \textbf{touch-vertices} of \( W' \).
\end{defi}

%Clearly,  when $Q_1=v_1$ and $Q_{k+1}=v_{k+1}$, an extended $R$-parity-switcher is precisely an $R$-parity-switcher.
The parity switcher plays a critical role in the  Hamilton-generated proofs as demonstrated in the following corollary.
For a path $P=v_1v_2\dots v_k$, denote by $v_iPv_j$ the subpath of $P$ that connects $v_i$ and $v_j$. 
For two internally disjoint paths $P$ and $Q$ with a common end $v$, denote by $PvQ$ the graph obtained by concatenating $P$ and $Q$ through their common end $v$.  
\begin{cor} \label{parity}
 Let $G$ be a Hamiltonian graph with an odd number $n$ of vertices, such that $\mathcal{C}_{n}(G) \ne \mathcal{C}(G)$.
    Let  $R$ be a subgraph described as in Lemma \ref{R-lem}. If $G$ has an $R$-parity switcher $$W=W(R, C_{(v_1, v_2, ... , v_{2k})}, P_1, P_2,\ldots, P_k, P_{k+1}, u_1, u_{k+1}),$$
then there does not exist a path $P$ connecting $u_1$ and $u_{k+1}$ with $\text{In}(P)=V(G) \setminus \bigcup_{i=1}^{k+1}V(P_i)$. 
\end{cor}
\begin{proof}
	Suppose to the contrary that there exists such a path $P$ with In$(P)=V(G) \setminus \bigcup_{i=1}^{k+1}V(P_i)$. Then we can construct two Hamilton cycles 
	from $P$ and $W$:	
	%$$\color{red}C_1=P_{k+1}+P+P_1+v_1v_2+P_2+v_{2k}v_{2k-1}+P_{3}+...+P_k+v_{k+2}v_{k+1}$$
 $$C_1=v_{k+1}P_{k+1}u_{k+1}Pu_1P_1v_1v_2P_2v_{2k}v_{2k-1}P_{3}v_{3}...v_kP_kv_{k+2}v_{k+1}$$
	and
	$$C_2=v_{k+1}P_{k+1}u_{k+1}Pu_1P_1v_1v_{2k}P_2v_{2}v_{3}P_{3}v_{2k-1}...v_{k+1}P_kv_kv_{k+1}.$$
It is clear that $E(C_1)\Delta E(C_2) = E(C)$. Note that $C$ has an odd number of edges in $R$. This implies that one of $C_1$ and $C_2$ has an odd number of edges in $R$, which contradicts the conclusion (b) of Lemma~\ref{R-lem}.
\end{proof}

Christoph, Nenadov, and Petrova~\cite{Nenadov} have established a key lemma that enables the efficient detection of short cycles with specific properties  in graphs. The lemma is stated as follows:
We define an $x$-$y$-path as a path connecting vertices $x$ and $y$. 
\begin{lem} [Cycle Lemma, \cite{Nenadov}] \label{Switcher}
    Let $G$ be a graph. Suppose subgraph $R \subseteq G$ satisfies the following conditions for some $ \ell \in \mathbf{N}:$

 (1) For every subset $S \subseteq V(G)$ of size $|S| \le 2\ell$ and any pair of vertices $x, y \in V(G)\setminus S$, there exists  an $x$-$y$-path in $R-S$ of length at most $\ell-1$;
 
(2) $ R \ne G$ and $R \ne G[A, B]$ for any partition $V (G) = A \cup B$.\\
 Then there exists an even cycle $C \subseteq G$ of length $|C|\le 2 \ell$ that contains an odd number of edges from $R$.

\end{lem}

As a slightly stronger version of a lemma of Krivelevich, Lee, and Sudakov~\cite{Sudakov}, Lemma \ref{Class} provides a classification of graphs with a minimum degree at least $(\frac 12 - o(1))n$.

\begin{lem} \label{Class}
 Let $\frac 1n \ll \beta \ll \alpha \ll 1$. Then for every graph $G$ on $n$ vertices with minimum degree at least $\left(\frac 12 - \beta\right) n$, one of the following holds:
 
 \begin{itemize}
     \item[(1)] $e(A,B) \ge \frac{\alpha}{2}n^2$ for all $A$ and $B$ with at least $\frac{1-\alpha}{2}n$ vertices (not necessarily disjoint).
     %$e(A, B) \ge  \alpha n^2$ for all A and B with at least $[\frac{n}{2}]$ vertices (not necessarily disjoint). {\color{blue} And this implies that : $e(A,B) \ge \frac{\alpha}{2}n^2$ for any all $A$ and $B$ with at least $\frac{1-\alpha}{2}n$ vertices (not necessarily disjoint).}
     \item[(2)] There exists a set A of size $\left(\frac 12 - 21 \alpha\right)n \le |A| \le \left(\frac 12 + 21 \alpha\right)n$ such that $e(A,\overline{A}) \le 4\alpha n^2$, and $\min\{\delta(G[A]), \delta(G[\bar{A}])\}\ge \frac n5$.     
     %the induced subgraphs on both $A$ and $\overline{A}$ has minimum degree at least $\frac n5$.
     \item[(3)] There exists a set A of size $\left(\frac 12 - 25 \alpha\right)n \le |A| \le \left(\frac 12 + 25 \alpha\right)n$  such that 
     $e(A, \bar{A})\ge (\frac 14 - 5\alpha)n^2$, $e(A) \le 6 \alpha n^2$, and $\delta(G[A, \bar{A}])\ge \frac n5$.
     %the bipartite graph induced by the edges between A and $\overline{A}$ has at least $(\frac 14 - 5\alpha)n^2$ edges and minimum degree $\frac n5$, and $e(A) \le 6 \alpha n^2$.
 \end{itemize}
\end{lem}

\begin{proof}
    Suppose that (1) does not hold. Then there exist two subsets $A_0$ and $B_0$ in $V(G)$ such that both of them have a size of at least $\frac{1-\alpha}{2}n$, and $e(A_0,B_0) < \frac{\alpha}2  n^2$. For convenience in the computation in the following proof, we may assume that $|A_0|=|B_0|=\left\lceil\frac {n}{2} \right\rceil$ and hence $e(A_0,B_0)\le \frac{\alpha}2  n^2+\frac{\alpha}2  n\cdot\frac n2< \alpha n^2$.
    %, {\color{red}which can be deduced from previous condition.}
    %\left\lceil\frac {n}{2} \right\rceil$ and $e(A_0,B_0) < \alpha n^2$.
Thus
    \begin{eqnarray}
        \sum_{u \in A_0 \cap B_0} d(u) &\le& \sum_{u \in A_0 \cap B_0}\left(|N(u)\cap (A_0 \cup B_0)|+|N(u)\setminus (A_0 \cup B_0)|\right)\nonumber\\
        &\le  & 2\alpha n^2+|A_0 \cap B_0|^2 \label{EQ: 1},   
    \end{eqnarray}
where the last inequality holds since $|N(u)\setminus (A_0 \cup B_0)|\le |V(G)|-|A_0|-|B_0|+|A_0\cap B_0|\le |A_0\cap B_0|$.
 Since $\delta(G)\ge (\frac 12- \beta)n$,  we have
 \begin{equation}\label{EQ: 2}
 	\sum_{u \in A_0 \cap B_0} d(u) \ge \left(\frac 12- \beta\right)n|A_0 \cap B_0|.
 \end{equation}   
From (\ref{EQ: 1}) and (\ref{EQ: 2}), we obtain that $|A_0 \cap B_0| \ge (\frac 12-5\alpha)n$ or $|A_0 \cap B_0| \le 5\alpha n$ (this can be guaranteed by taking $\beta<\frac 12\alpha<\frac1{350}$).

If $|A_0 \cap B_0| \le 5\alpha n$, let $X' = A_0 \setminus B_0$ and $Y' = \overline{X'}=B_0\cup\overline{A_0\cup B_0}$, then $(\frac 12-5\alpha)n \le |X'|,|Y'| \le (\frac 12+5\alpha)n$. 
Note that $|\overline{A_0 \cup B_0}|=|V(G)|-|A_0|-|B_0|+|A_0\cap B_0|\le |A_0\cap B_0|\le 5\alpha n$. Thus 
$$e(X',Y') = e(A_0 \setminus B_0, B_0) + e(A_0 \setminus B_0, \overline{A_0 \cup B_0}) < \alpha n^2 + 5\alpha n|X'| \le 4\alpha n^2,$$
and  $$e(X')=\frac 12\left(\sum_{u\in X'}d(u)-e(X', Y')\right) \ge \left(\frac 14-\frac{\beta}2\right) n|X'| - 2\alpha n^2 \ge \left(\frac 18-4\alpha\right) n^2,$$
and similarly, $e(Y') \ge (\frac 18-4\alpha) n^2$.
Let $W=\{x\in X' : |N(x)\cap X'|\le \frac{6}{25}n\}$ and $Z=\{y\in Y' : |N(y)\cap Y'|\le \frac{6}{25}n\}$.    
%    Furthermore, we can let the set of vertices in $X'$ such that has at most $0.24 n$ neighbors in $X'$ be $W$ and set of vertices in $Y'$ such that has at most $0.24 n$ neighbors in $Y'$ be $Z$, then by
Then $|N(x)\cap Y'|\ge (\frac 12-\beta)n-\frac 6{25}n>\frac 14n$ for each $x\in W$,  and similarly, $|N(y)\cap X'|\ge \frac 14n$ for any $y\in Z$.
Thus $$\frac 14n|W| \le e(W,Y') \le e(X',Y') \le 4\alpha n^2,$$ and similarly, $\frac 14n|Z| \le e(Z,X') \le e(X',Y') \le 4\alpha n^2$.
 Therefore,   we have
    $|W|\le 16 \alpha n$ and $|Z| \le 16 \alpha n$.
    Let $X=(X'\setminus W)\cup Z$ and $Y=(Y' \setminus Z)\cup W$. Then 
        $$\left(\frac 12 - 21\alpha\right) n \le |X|,|Y| \le \left(\frac 12+ 21\alpha\right) n,$$
and $$\delta(G[X])\ge \frac{6n}{25}-16\alpha n\ge \frac n5$$
and similarly, $\delta(G[Y]) \ge \frac{n}{5}$,
%{\color{red}$$ \delta(G[X]), \delta(G[Y]) \ge 0.24n \ge \frac{n}{5},???$$}
    and
    $$e(X,Y) \le e(X',Y') \le 4\alpha n^2.$$
Therefore, the statement (2) holds.

     If $|A_0 \cap B_0| \ge (\frac 12-5\alpha)n$, let $X' = A_0 \cap B_0$ and $Y' = \overline{A_0 \cap B_0}$, then we have $(\frac 12-5\alpha)n \le |X'|,|Y'| \le (\frac 12+5\alpha)n$, and $e(X') \le e(A_0,B_0) \le \alpha n^2$, and $$e(X',Y') \ge \left(\frac 12-\beta\right)n|X'|-2\alpha n^2 \ge \left(\frac 14-5\alpha\right)n^2.$$
Let $W=\{x\in X' : |N(x)\cap Y'|\le \frac{6n}{25}\}$ and $Z=\{y\in Y' : |N(y)\cap X'|\le \frac{6n}{25}\}$.     
%     Furthermore, we can let the set of vertices in $X'$ such that has at most $0.24 n$ neighbors in $Y'$ be $W$ and set of vertices in $Y'$ such that has at most $0.24 n$ neighbors in $X'$ be $Z$, 
Then by
    $$
      e_{G^c}(W, Y') \le |X'||Y'|-e(X',Y') \le \frac{n^2}{4}-\left(\frac{1}{4}-5\alpha\right)n^2=5\alpha n^2
     $$
     and
     $$e_{G^c}(W, Y')\ge \left(|Y'|-\frac {6n}{25}\right)|W|\ge \frac n4 |W|,$$
     we have
     $|W|\le 20 \alpha n$. Similarly, we have $|Z|\le 20 \alpha n$.
     Let $X=(X'\setminus W)\cup Z$ and $Y=(Y'\setminus Z)\cup W$. 
     Then 
     $$\left(\frac 12 - 25\alpha\right) n \le |X|,|Y| \le \left(\frac 12+ 25\alpha\right) n,$$ 
     and  
     \begin{align*}
         &\delta_X(Y'\setminus Z) \ge \frac{6n}{25}-|Z| \ge \frac{6n}{25}-20\alpha n\ge\frac n5,\\
     &\delta_X(W)  \ge \delta(G)-\frac{6n}{25}-20\alpha n \ge 
     (\frac{13}{50}-20\alpha -\beta)n \ge \frac n5.
     \end{align*}
     Hence $\delta_X(Y)=\min\{\delta_X(Y'\setminus Z),\delta_X(W)\} \ge \frac n5$, and similarly, $\delta_Y(X)\ge \frac n5$.
%     $$
%     \delta_Y(X) \ge 0.24n \ge \frac{n}{5},$$  $$\delta_X(Y), \delta_Y(X) \ge 0.24n \ge \frac{n}{5},$$ 
Note that  $|N(z)\cap X|\le\frac{6n}{25}+20\alpha n\le\frac n4$ for $z\in Z$. We have $e(X) \le e(X') + \frac{n}{4} |Z| \le 6\alpha n^2$. 
Clearly,      $e(X,Y) \ge e(X',Y') \ge (1/4-5\alpha)n^2.$ Therefore, the statement (3) holds.

\end{proof}

We call a graph $G$ is \textit{$(m,d)$-connected} if for any $U\subseteq V(G)$ of size at most $m$, the resulting graph remains connected and has a diameter of at most $d$ after deleting $U$. It is easy to check that for any $(m,d)$-connected $G$ and $U_0 \subseteq V(G)$ of size $m'<m$, $G-U_0$ is $(m-m',d)$-connected.
 
 \begin{lem} \label{bidiameter}
    Let $m,d \in \mathbb{N}^*$ and $\ell \le \lceil\frac{m}{d+1}\rceil$. If $G$ is $(m,d)$-connected, for any $\ell$ disjoint pairs of vertices, $(u_1,v_1),(u_2,v_2),\dots,(u_{\ell},v_{\ell})$, in $G$, we can find $\ell$ disjoint paths $P_1,P_2,\dots,P_{\ell}$ of length at most $d$ such that $P_i$ connects $u_i$ and $v_i$ for $1\le i\le\ell$.

  \end{lem}
  \begin{proof}
    The desired $\ell$ disjoint paths can be constructed greedily one by one. 
      Assume that we have constructed $j$ desired disjoint paths $P_1,P_2,\dots,P_{j}$ with $|E(P_i)| \le d$ for all $i \le j\le \ell-1$.
    Let $U=(\cup_{1=i}^jV(P_i))\bigcup(\cup_{k=j+2}^\ell\{u_k, v_k\})$. Then $|U| \le (d+1)j+2(\ell-j-1)< (d+1)\ell \le m$. Hence $G-U$ has a diameter of at most $d$. This implies that there exists a path $P_{j+1}$ of length at most $d$ such that $P_{j+1}$ connects $u_{j+1}$ and $v_{j+1}$, and $V(P_{j+1}) \subseteq V(G)\setminus U$. Therefore, we can find $\ell$  such paths $P_1,P_2,\dots,P_{\ell}$ as desired. 
\end{proof} 

The following two lemmas provide minimum degree conditions to guarantee a Hamilton cycle passing through a set of specific edges. 
\begin{lem}[Posa, \cite{posa1963}]\label{dirac:HC}
    Let $\mathcal{F}$ be a UDP of $k > 1$ edges in a graph $G$ of order $n$. If $\delta(G) \ge \frac{n+k}{2}$ then $G$ has a Hamilton cycle containing $\mathcal{F}$.
\end{lem}
For a bipartite graph $G$ with partite sets $X$ and $Y$, let $\sigma_{1,1}(G) = \min \{ d_G(x) + d_G(y) : x \in X, y \in Y, xy \notin E(G) \}$
if \( G[X, Y] \) is not complete; otherwise, \( \sigma_{1,1}(G) = \infty \).
For a positive integer $k$, a {\em $k$-matching} in a graph is a matching of size $k$.
For a given $k$-matching $M$  in a graph, we call a Hamilton cycle of $G$ as a \textit{Hamilton $M$-cycle} if it passes through all edges of $M$.
\begin{thm}[\cite{Fuji24}]
    Let \( G \) be a balanced bipartite graph of order \( 2n \), and \( M \) be a \( k \)-matching in \( G \), where \( k \geq 1 \). Then \( G \) contains a Hamilton \( M \)-cycle if
\[
\sigma_{1,1}(G) \geq 
\begin{cases} 
n+2 & \text{if } k \in \{n, n-1\} \text{ or } 1 \leq k \leq 4, \\ 
2n - k & \text{if } \frac{2n}{3} \leq k \leq n-2, \\ 
n + \frac{k}{2} & \text{if } 5 \leq k \leq \frac{2n-1}{3}.
\end{cases}
\]
Moreover, the lower bound is sharp for every \( k \) with \( 1 \leq k \leq n \).
\end{thm}

From the preceding two lemmas, we establish  a key technical lemma  that assists us in identifying a Hamilton cycle traversing  all specified paths in graphs satisfying certain favorable degree conditions, even when these paths contain a  sublinear number of edges.
\begin{lem} \label{Fixpair}
    Suppose that $\frac{1}{n}\ll \varepsilon_1$, $\varepsilon_2 \ll 1$ and $\ell \le \varepsilon_2 n$.
    
    (1) If graph $G$ has $n$ vertices and $\delta(G) \ge (1-\varepsilon_1) n$, then for any $\ell$ disjoint pairs of vertices $\{(u_1,v_1),(u_2,v_2),\dots,(u_\ell,v_\ell)\}$, there exists a Hamilton cycle in 
    \begin{equation*}
        G + \bigcup_{1 \le i \le \ell}\{u_iv_i\}
    \end{equation*}    
    that contains all the above $\ell$ edges.

    (2) If $G$ is a bipartite graph with $n$ vertices, partitioned into two sets $X$ and $Y$, such that: 
    $ (\frac 12-\varepsilon_1)  n \le |X| \le  (\frac 12+\varepsilon_1) n$ and 
    $ \delta(G) \ge (\frac{1}{2}-\varepsilon_1)n$, then for any  $\ell$ disjoint pairs of vertices $\{(u_1,v_1),(u_2,v_2),\dots,(u_\ell,v_\ell)\}$ satisfying     
    \begin{equation}\label{balanced}
        \left|X\right|- \frac{1}{2}\left|\bigcup_{1 \le i \le \ell} \{u_i,v_i\} \cap X\right|=|Y|- \frac{1}{2}\left|\bigcup_{1 \le i \le \ell} \{u_i,v_i\} \cap Y\right|,
    \end{equation}
 there exists a Hamilton cycle in 
    \begin{equation*}
        G +\bigcup_{1 \le i \le \ell}\{u_iv_i\}
    \end{equation*}
containing all the above $\ell$ edges.
\end{lem}
\begin{proof}
(1) Note that $\delta(G) \geq (1-\varepsilon_1)n \ge \frac{n+\ell}{2}$.
Then we can prove this case by applying Theorem~\ref{dirac:HC} for $G + \bigcup_{1 \leq i \leq \ell} \{u_i v_i\}$ with $\mathcal{F}$ consisting of all $\ell$ edges in $\bigcup_{1 \leq i \leq \ell} \{u_i v_i\}$.

(2) Since \( \left(\frac{1}{2} - \varepsilon_1\right) n \leq |X|, |Y| \leq \left(\frac{1}{2} + \varepsilon_1\right) n \) and \( \delta(G) \geq \left(\frac{1}{2} - \varepsilon_1\right) n \),    
for any two vertices in the same part, they have at least
\[
2\delta(G) - \max\{|X|, |Y|\} \geq \left(\frac{1}{2} - 3\varepsilon_1\right)n
\]
common neighbors in the other part.
For any two vertices \( u \) and \( v \) in different parts, we choose a neighbor of \( u \), denoted as \( u' \), which is in the same part as \( v \). Then \( u' \) and \( v \) have at least \( \left(\frac{1}{2} - 3\varepsilon_1\right)n \) common neighbors. 
Choose $\varepsilon_3$ such that $\frac1n \ll \varepsilon_1 \ll \varepsilon_2 \ll \varepsilon_3< \frac{1}{2} - 3\varepsilon_1$.
Hence, \( G \) is \((\varepsilon_3 n, 3)\)-connected.

Choose an additional vertex pair $(v_0, u_{\ell+1})$ with \( v_0 \in X \) and \( u_{\ell +1} \in Y \) distinct from all pairs \( (u_i, v_i) \)'s.  
Given that \( \ell \leq \varepsilon_2 n \leq \frac{\varepsilon_3 n}{3} \), we can apply Lemma~\ref{bidiameter} to find \( \ell + 1 \) disjoint paths \( P_0, P_1, P_2, \dots, P_{\ell} \), each of length at most \( 3 \), such that \( P_i \) connects \( v_i \) and \( u_{i+1} \) for \( 0 \leq i \leq \ell \).

Since all paths $P_i$ are in bipartite graph $G$ for any \( i \in \{0,1,2,\dots,\ell\} \), we have:
\begin{itemize}
    \item \( |V(P_i) \cap X| - |V(P_i) \cap Y| = 1 \) if the endpoints of $P_i$ are both in \( X \);
    \item \( |V(P_i) \cap X| - |V(P_i) \cap Y| = 0 \) if the endpoints of $P_i$ are in different vertex parts;
    \item \( |V(P_i) \cap X| - |V(P_i) \cap Y| = -1 \) if the endpoints of $P_i$ are both in \( Y \).
\end{itemize}
In all, we have 
\begin{equation}\label{balanced:2}
        \left|V(P_i) \cap X\right|- \frac{1}{2}\left| \{v_i,u_{i+1}\} \cap X\right|=\left|V(P_i) \cap Y\right|- \frac{1}{2}\left| \{v_i,u_{i+1}\} \cap Y\right|.
    \end{equation}
By concatenating these paths along with the edges $\{u_i v_i\}$ for $1 \leq i \leq \ell$, we obtain a path
$
P_{v_0u_{\ell +1}} = v_0P_0u_1 v_1 P_1 u_2 \dots v_{\ell-1} P_{\ell-1} u_{\ell} v_{\ell}P_{\ell}u_{\ell +1}
$
in graph $G + \bigcup_{1 \leq i \leq \ell} \{u_i v_i\}$. This path contains at most $4(\ell+2) \leq 5\varepsilon_2 n$
vertices.

Let $G'=G- \text{In}(P_{v_0u_{\ell +1}})+\{v_0u_{\ell +1}\}$. Then $|V(G')|=n'=n-|\text{In}(P_{v_0u_{\ell +1}})| \ge (1-5\varepsilon_2)n$ and it is even. This is because, according to (\ref{balanced}) and (\ref{balanced:2}), we have $|X\cap V(G')|=|X\setminus \text{In}(P_{v_0u_{\ell +1}})|=|Y\setminus \text{In}(P_{v_0u_{\ell +1}})|=|Y\cap V(G')|$.
Note that
$\delta(G') \ge \delta(G)-|\text{In}(P_{v_0u_{\ell +1}})| > (\frac12-6\varepsilon_2)n'$. Hence $\sigma_{1,1}(G')\ge 2\delta(G') \ge (1-12\varepsilon_2)n' \ge \frac{n'}{2}+2$.
Apply Theorem~\ref{dirac:HC} to $G'$ with $1$-matching $M=\{v_0u_{\ell +1}\}$, we obtain a Hamiltonian cycle in $G'$ passing through the edge $v_0u_{\ell +1}$. Thus  we have a Hamiltonian path $P'$ in $G'$ connects $v_0$ and $u_{\ell +1}$.
Therefore, the cycle $C' = v_0 P_{v_0 u_{\ell+1}} u_{\ell+1} P' v_0$ is a Hamilton cycle as we desired.

\end{proof}  

Our main contribution consists of the following three lemmas, each handling a distinct scenario based on the classification established in Lemma~\ref{Class}. For a positive integer $k$, a {\em $k$-cycle} in a graph is a cycle of length $k$.
\begin{lem} \label{biclique} 
    Suppose $n$ is an odd number and $\frac{1}{n}\ll \alpha \ll \eta$, and let $G$ be a graph on $n$ vertices with a partition $V(G) = X \cup Y \cup Z$ such that:
    \begin{equation*}
        \left(\frac{1}{2}- \alpha\right)n \le |X|, |Y| \le \left(\frac{1}{2}+ \alpha\right)n, \, |Z| \le \alpha n,
    \end{equation*}
 $e(X,Y) \le \alpha n^2$, and $ \max\{\Delta_X(Y),\Delta_Y(X)\} \le \eta n$,
    \begin{equation*}
       \min\{\delta(Y), \delta(X) \}\ge \left(\frac{1}{2}- 2\eta\right)n, \,\,     \min\{\delta_Y(Z),    \delta_{X}(Z)\} \ge \frac34 \eta n,
    \end{equation*}
 then the following holds:
If 
   \begin{eqnarray} \label{cod:biclique}
       \frac{4}{3}|Z|+m(X,Y) \ge \frac{10}{3},
   \end{eqnarray}
     then
   $\mathcal{C}_{n}(G)=\mathcal{C}(G)$. Here, $m(X,Y)$ is the matching number of the bipartite graph $G[X,Y]$.
\end{lem}

\noindent\textbf{Remark.} We emphasize that the satisfaction of the inequality (\ref{cod:biclique}) is necessary. This requirement is demonstrated by Construction~A: while both graphs $G_1$ and $G_3$ fails to be Hamilton-generated as demonstrated in Proposition~\ref{example}. They satisfy all criteria specified in Lemma~\ref{biclique} exception for the inequality (\ref{cod:biclique}) (in fact, the corresponding values of these graphs are $\frac{8}{3}$ and $3$, respectively).

\begin{proof}
    Suppose that $\mathcal{C}(G) \neq \mathcal{C}_n(G)$. According to Lemma \ref{R-lem}, there exists a subgraph $R$ of $G$ satisfying 
    \begin{itemize}
       \item[(a)] $R \ne G$,
       \item[(b)] every Hamilton cycle in $G$ contains an even number of edges from $R$, and
       \item[(c)] for every partition $V (G) = A \cup B$, we have $e_{R}(A,B) \ge \frac 12e_{G}(A,B)$ and $R \ne G[A,B]$.
   \end{itemize}
From (c), we have $\delta(R) \ge \frac{1}{2}\delta(G)$. 
   Choose $1/n \ll \alpha \ll \eta \ll \sigma \ll 1$. Denote $Z$ as $\{z_1,z_2,...,z_{|Z|}\}$. 
 As both $\delta_Y(Z)$ and $\delta_{X}(Z)$ are at least $\frac34 \eta n(\ge 3|Z|)$,  by the pigeonhole principle, for each $z_i$, we
  can  choose $x_i^1,x_i^2 \in N_G(z_i) \cap X$ and $y_i^1,y_i^2 \in N_G(z_i) \cap Y$ such that either $x_i^1 z_i, x_i^2z_i\in E(R)$ or $x_i^1 z_i, x_i^2z_i\in E(G) \setminus E(R)$, and either $y_i^1 z_i, y_i^2 z_i\in E(R)$ or $y_i^1 z_i, y_i^2 z_i\in E(G) \setminus E(R)$. 
  %{\color{brown}Furthermore, for any set $S$ with at most $\frac{\eta n}{2}$ vertices, we can assume that all the above chosen  neighbors are disjoint from $S$ due to the same reason that $\delta_{Y\setminus S}(Z), \delta_{X\setminus S}(Z) \ge \frac34 \eta n - \frac{1}{2}\eta n = \frac{1}{4}\eta n>3|Z|$.}

\begin{claim} \label{connected}
Both $G[X]$ and $G[Y]$ are $(\frac{1}{2}\eta n, 2)$-connected. Moreover,  $G[X\cup Z]$ and $G[Y\cup Z]$ are each $(\frac{1}{2}\eta n,4)$-connected. 
      %{\color{brown} Directly, $G[X]$ and $G[Y]$ are $((\frac{\eta }{2}-\alpha) n,4)$-connected.}
      %For any subset $U$ in $G[X \cup Z]$  (resp. in $G[Y \cup Z]$) of size at most $\frac{\eta n}{2}$, the resulting graph remains connected and has a diameter of at most $4$ after deleting $U$.
  \end{claim}
  \begin{proof}
    We only need to prove the claim for $G[X]$ and $G[X\cup Z]$ by the symmetry of $X$ and $Y$. Choose a subset $U\subseteq X\cup Z$ such that $|U| \le \frac{\eta n}{2}$ and $u,v \in (X \cup Z) \setminus U$ randomly. Since $\delta(X)\ge \left(\frac{1}{2}- 2\eta\right)n$ and $\delta_{X}(Z)\ge \frac{3}{4}\eta n$, we have 
\begin{equation}\label{EQ: e1-mindeg}
    \delta(X\setminus U) \ge \delta(X)-|U| \ge \left(\frac{1}{2}-2\eta\right)n-\frac{\eta n}{2}  \ge \frac{1/2-3\eta}{1/2+\alpha}|X| \ge (1-7\eta)|X|
\end{equation}
and 
\begin{equation}\label{EQ: e2-mindegZ}
  \delta_{X\setminus U}(Z) \ge \delta_X(Z)-|U| \ge \frac{\eta n}{4}.  
\end{equation}
%$$\delta_{X\setminus U}(Z) \ge \delta_X(Z)-|U| \ge \frac{\eta n}{4}.$$
Inequality (\ref{EQ: e1-mindeg})  implies that if  $u$ and $v$ are both in $X\setminus U$, then they have a common neighbor in $X\setminus U$, i.e. there exists a path of length at most $2$ connecting $u$ and $v$. 
Inequality (\ref{EQ: e1-mindeg}) together (\ref{EQ: e2-mindegZ}) imply that  we always can choose two neighbors of $u$ and $v$ in $X\setminus U$, respectively, and these two neighbors have a common neighbor, i.e. there exists a path of length at most $4$ connecting $u$ and $v$. As $u$ and $v$ are chosen randomly, $G[X]$ is $(\frac{1}{2}\eta n, 2)$-connected and $G[(X\cup Z)\setminus U]$ is $(\frac{1}{2}\eta n, 4)$-connected.
\end{proof}

From the preceding construction of $u$-$v$-path, we have the following claim.
\begin{claim}\label{CL: differlengthpath}
Let $u$ and $v$ be vertices in $X\cup Z$ or in $Y\cup Z$. Then there exists a  $u$-$v$-path of length $\ell$ for $\ell\in\{2,3,4\}$ if $\{u,v\}$ is a pair from $X$ or $Y$, of length $\ell$ for $\ell\in\{3,4\}$ if one vertex in $X\cup Y$ and the other in $Z$, and of length four if $\{u, v\}$ is a pair from $Z$.
\end{claim}

\begin{claim} \label{noR-parity}
There do not exist an \( R \)-parity-switcher \( W \) and an edge \( e \) such that \( W \) has at most \( \frac{5}{12}\eta n + 100 \) vertices, with two touch-vertices \(u_1 \in X \cup Z \) and \( u_{k+1} \in Y \cup Z \), and \( e \in E_G(X \cup Z, Y \cup Z) \) is not incident with \( V(W) \).
       
\end{claim}

    \begin{proof}
           Suppose to the contrary that there exists such an $R$-parity-switcher $$W=W(R,C_{(v_1, v_2, \ldots, v_{2k})}, P_1,P_2\dots,P_k,P_{k+1}, u_1, u_{k+1})$$ and an edge $e=xy$ with $x \in (X\cup Z)\setminus V(W)$ and $y \in (Y\cup Z)\setminus V(W)$.
            %, such that $W=W(R, C_{(v_1, v_2, ... , v_{2k})}, P_2,\ldots, P_k, v_1, v_{k+1})$ is an $R$-parity-switcher, and $Q_i$ connects $v_i$ and $v_i'$ for $i \in \{1,k+1\}$. 
            Let $U=V(W) \cup \{x,y\}$ and let $X'=X\setminus U$, $Y'=Y\setminus U$, and $Z'=Z\setminus U$. Then 
           $|U| \leq |V(W)| + 2 \leq \frac{5}{12}\eta n + 102$. Hence
      \begin{align*}
   \min\{ \delta_{X'}(X' \cup Z'), \delta_{Y'}(Y' \cup Z') \}&\geq \frac{3}{4} \eta n - \left(\frac{5}{12} \eta n + 102\right) \geq \frac{1}{4} \eta n.
      \end{align*}
      This implies that we can
      choose $x'$ and $w_1$ as neighbors of $x$ and $u_1$ in $X'$, respectively, $y'$ and $w_{k+1}$ as neighbors of $y$ and $u_{k+1}$ in $Y'$, respectively.
      %$x' \in (N_G(x) \cap X)\setminus V(W')$.  
      %Let $X',Y',Z'$ be the remaining set of $X,Y,Z$ by deleting $(V(W)\cup V(Q_1)\cup V(Q_{k+1}))\setminus \{v_1',v_{k+1}'\}$ respectively.
      
      If there exist two paths $P_{w_1x'}$ between $w_1$ and $x'$ with vertex set $X'$,  and $P_{w_{k+1}y}$ between $w_{k+1}$ and $y$ with vertex set $Y'\cup Z'$, then we can construct a path 
      $$P_{v_1v_{k+1}}=v_1P_1u_1w_1P_{w_1x'}x'xyy'P_{w_{k+1}y'}w_{k+1}u_{k+1}P_{k+1}v_{k+1}$$ such that $V(P_{v_1v_{k+1}})=V(G)\setminus \bigcup_{i=2}^k V(P_i)$. This construction contradicts Corollary~\ref{parity}. %So we only need to find such two paths. 
      Thus, we only need to construct such two paths \( P_{w_1x'} \) and \( P_{w_{k+1}y'} \). Indeed, we only need to construct \( P_{w_{k+1}y'} \). This is because \( X' \) and \( Y' \cup Z' \) are symmetric in a certain sense, even though \( Z' \) is an additional part.
      
      %Now we are going to construct such $P_{w_1x'}$ and $P_{w_{k+1}y'}$.
      %{\color{red}In sure, we only need to construct $P_{v_{k+1}'y}$ since $X'$ and $Y'\cup Z'$ are symmetric in some sense, although there is an additional part $Z'$.}
      
      Recall that for any $i \in \{1,2,\dots,|Z|\}$, we have chosen two neighbors $y_i^1,y_i^2$ of $z_i$  in $Y$. We may additionally assume that
      $Z'=\{z_1,z_2,\dots,z_{|Z'|}\}$
      and any $y_i^j$ as above locates at $Y'\setminus \{v_{k+1}',y\}$.
      This is because for any $z \in Z$, we have 
      \begin{align*}
          |(N_G(z) \cap Y')\setminus\{y',w_{k+1}\}| \ge \delta_{Y'}(Y'\cup Z')-2 \ge \frac14\eta n - 2 \ge 3\alpha n  \ge 3 |Z'|.
      \end{align*}
 Note that $\delta(Y') \ge \delta(Y)-|U| \ge (\frac12-2\eta)n-\frac{5}{12}\eta n-102 \ge (1-6\eta)|Y'|$.  Apply Lemma \ref{Fixpair} (1) to
        $G[Y']$ with $\varepsilon_1=6\eta$ and $\varepsilon_2=2\alpha$, and pairs of vertices 
        $$
        \{(y',w_{k+1})\}\cup \left(\cup_{i=1}^{|Z'|} \{(y_i^1,y_i^2)\}\right).
        $$
Clearly, the number of pairs $|Z'|+1 \le  |Z|+1 \le 2\alpha n$.
 Then we have a Hamilton cycle $C$ containing all edges $y'w_{k+1}$ and $y_i^1y_i^2$ for $1\le i\le |Z'|$ in  $G[Y]+\{y'w_{k+1}, y_i^1y_i^2, 1\le i\le |Z'|\}$.  Therefore, we obtain a path $P_{y'w_{k+1}}$  that covers exactly all vertices in $Y'\cup Z'$ from $C$ by replacing each $y_i^1y_i^2$ by $y_i^1z_iy_i^2$ and deleting the edge $y'w_{k+1}$.
%. \footnote{where $P_{y'w_{k+1}}$ can be obtained from the resulting Hamilton cycle containing $\{y'w_{k+1}, y_i^1y_i^2, 1\le i\le |Z'|\}$ in  $G[Y]+\{y'w_{k+1}, y_i^1y_i^2, 1\le i\le |Z'|\}$ by replacing each $y_i^1y_i^2$ by $y_i^1z_iy_i^2$ and deleting the edge $y'w_{k+1}$.}

\end{proof}

Let $M_0$ be a matching in $G[X \cup Z, Y]$.
    
\begin{claim} \label{even-cycle}
       If $|M_0|=2$, then there does not exist an even cycle in $G[X \cup Z]$ of length at most $\frac{1}{6}\eta n$ such that it contains at most one vertex in $V(M_0)$ and has an odd number of edges in $E(R)$.
\end{claim}
\begin{proof}    
Suppose not, let $C$ be such a cycle with length $2h \le \frac{1}{6}\eta n$ vertices.
      Then $G[(X\cup Z)\setminus V(C)]$ is connected by Claim~\ref{connected}. Since $|V(C) \cap V(M_0)| \le 1$, there is at least one edge between $(X \cup Z)\setminus V(C)$ and $Y$. 
      %This insures that $G-V(C)$ remains connected.
      Choose vertices $a \in V(C)$  and $b \in V(M_0)\cap (X\cup Z)$ (if $|V(M_0)\cap V(C)|=1$ then choose $a=b$). 
      Let $U=(V(C) \cup V(M_0))\setminus \{a,b\}$. Then $|U\cap (X\cup Z)|\le \frac 16\eta n+2 \le \frac12\eta n $. 
      By Claim~\ref{connected}, $G[X\cup Z]$ is $(\frac{1}{2}\eta n,4)$-connected. Hence 
      we have a path $P_{ab}$ of length at most 4 connecting $a$ and $b$ in $G[(X\cup Z)\setminus U]$. Clearly, $V(P_{ab})\cap V(M_0)=\{b\}$.
      Denote $M_0$ as $\{bb', xy\}$,  where $x \in X \cup Z$ and $y \in Y$.
%      \iffalse
%      Since $\delta_Y(Z), \delta_{X}(Z) \ge \frac34 \eta n,$ we may further assume that
%      $$\{a,b,c,d,b'\} \cap (\cup_{i=1}^{|Z\setminus V(C)|} \{x_i^1,x_i^2,y_i^1,y_i^2\})=\emptyset.$$
%Note that 
%\begin{equation}\label{mindelta_Y}
%    \delta(Y)\ge \left(\frac 12-2\eta\right)n \ge \frac{1/2-2\eta}{1/2+\alpha}|Y| \ge (1-5\eta)|Y|.
%\end{equation}
%Applying Lemma \ref{Fixpair} (1) to
%        $G[Y]$ with $\varepsilon_1=5\eta$ and $\varepsilon_2=2\alpha$, and pairs of vertices 
%        $$
%        \{(b',d)\}\cup \left(\cup_{i=1}^{\ell} \{(y_i^1,y_i^2)\}\right),
%        $$
% where the number of pairs $\ell+1=|Z \setminus (V(C)\cup \{c\} \cup V(P_1))|+1 \le  |Z|+1 \le 2\alpha n$,
%  we obtain a path $P_{b'd}$ from $b'$ to $d$ that covers exactly all vertices in $(Y \cup Z)\setminus (V(C)\cup \{c\} \cup V(P_{ab}))$ (where $P_{b'd}$ can be obtained from the resulting Hamilton cycle  containing $\{b'd, y_i^1y_i^2, 1\le i\le \ell\}$ in  $G[Y]+\{b'd, y_i^1y_i^2, 1\le i\le \ell\}$ by replacing each $y_i^1y_i^2$ by $y_i^1z_iy_i^2$ and deleting the edge $b'd$). 
%  \fi
Label $C$ as $(v_1,v_2,...,v_{h},v_{h+1},...,v_{2h},v_1)$, where $a=v_1$. Let $P_1=P_{ab}+bb'$ and $P_{h+1}=\{v_{h+1}\}$. Then $P_1$ has a length at most $5$.
Since $G[X\cup Z]$ is $(\frac12\eta n,4)$-connected and $|Z \cup V(P_{ab})\cup \{v_{h+1},x\}| \le \alpha n + 7 < \frac{1}{12}\eta n$, we have $G[(X\cup V(C))\setminus (V(P_{ab})\cup \{v_{h+1},x\})]$ is $(\frac{5}{12}\eta n,4)$-connected by the definition of $(\cdot , \cdot)$-connected. 
By applying Lemma~\ref{bidiameter} to $G[(X\cup V(C))\setminus (V(P_{ab})\cup \{v_{h+1},x\})]$, we can find $h-1$ (where $h-1\le \frac{1}{12}\eta n -1 < \lceil{\frac{5\eta n}{12}}/{5}\rceil)$ disjoint paths $P_{2},P_{3},\dots,P_{h}$ such that $P_i$ connects $v_i$ and $v_{2h+2-i}$ and $|E(P_i)| \le 4$ for $2 \le i \le h$. Thus we have  an $R$-parity-switcher $W=W(R,C, P_1, P_2,\ldots, P_h, P_{h+1},  b', v_{h+1})$. 
Note that $|V(W)| \le 6+1+5(h-1)=5h+2 \le \frac{5}{12}\eta n+2 \le \frac{5}{12}\eta n +100$ and $V(W)\cap \{x,y\}=\emptyset$. The $R$-parity switcher $W$ and the edge $xy$ together forms a contradiction to Claim~\ref{noR-parity}.

    \end{proof}

    \begin{claim} \label{claim2}
For any $U \in \{X,Y\}$, if $E(G[U]) \setminus E(R) \neq \emptyset$ and $|M_0|=3$, then there exist a subset  $A \subseteq U$ and a partition $Z=Z_A \cup Z_{U\setminus A}$ such that
\begin{itemize}
    \item $(\frac 12-25\sigma)|U| \le |A| \le (\frac 12+25\sigma)|U|$;
    \item $R[U \cup Z]=G[A\cup Z_A,(U\setminus A)\cup Z_{U\setminus A}]$;
    \item $R[A,U\setminus A]$ is $(\sqrt{5\sigma}n,4)$-connected.
\end{itemize}

        %{\color{blue} $R[X]$ forms a cut of $G[X]$ with two parts $A$ and $B$. Moreover, %$(\frac12-25\sigma)|X| \le |A|,|B| \le (\frac12+25\sigma)|X|$
        %$(\frac 12-25\sigma)|X| \le |A|,|B| \le (\frac 12+25\sigma)|X|$, 
       %$\delta(R[A, B])\ge \frac{|X|}5 \ge \frac n{11}$, $e_R(A,B) \ge (\frac{1}{4}-5\sigma)|X|^2 \ge (\frac{1}{16}-2\sigma)n^2$.}
    \end{claim}

    \begin{proof}
   We only need to prove the claim for $U=X$ by the symmetry of $X$ and $Y$.
    According to Claim \ref{even-cycle}, there is no even cycle $C$ in $G[X\cup Z]$ satisfying: (a) $|V(C)|\le \frac{\eta n}{6}$, and (b) $|V(C)\cap V(M')|\le 1$ for any 2-subset $M'$ of $M_0$ (which implies  $|V(C) \cap V(M_0)| \le 2$),  and (c) $|E(C) \cap E(R)|$ is odd.
   %     \begin{itemize}
   %         \item[(a)] $|V(C)|\le \eta n$;
    %        \item[(b)] {\color{red}$|V(C) \cap V(M_0)| \le 2$;}
     %       \item[(c)] $|E(C) \cap E(R)|$ is odd.
   %     \end{itemize}
     We will prove this claim in two steps.

\noindent
\textbf{Step I.} We show that there exists a subset $A \subseteq X$ such that $R[X]=G[A,X\setminus A]$ and $(\frac 12-25\sigma)|X| \le |A| \le (\frac 12+25\sigma)|X|$.

Recall that $|Z| \le \alpha n$, $\Delta_Y(X) \le \eta n$, and $\delta(R)\ge \frac12\delta(G)\ge\frac 12\delta(X)$. Then 
  $$\delta(R[X]) \ge\delta(R)-\Delta_Y(X)-|Z|\ge  \frac{1}{2} \delta(X) - 2\eta n \ge \left(\frac{1}{4}-3\eta\right)n \ge \left(\frac{1}{2}-7\eta\right)|X|,$$
        the last inequality holds since $|X| \le (1/2+\alpha)n$.
Hence we may apply Lemma~\ref{Class} to $R[X]$ for $7\eta\ll\sigma$. 

If for any $A,B \subseteq X$ with $|A|,|B|\ge \frac{1-\sigma}{2}|X|$, we have $e_R(A,B)\ge \frac{\sigma}{2} |X|^2$~(Lemma~\ref{Class} (1)).
        Choose an edge $e=ab \in E(G[X]) \setminus E(R)$. Then
        $$|N_R(a) \cap X|, |N_R(b) \cap X| \ge \delta(R[X])\ge \left(\frac 12-7\eta\right)|X|.$$
 Hence $|N_R(a) \cap (X \setminus V(M_0))|,|N_R(b) \cap (X \setminus V(M_0))| \ge (\frac12-7\eta)|X|-3 \ge \frac{1-\sigma}{2}|X|$
 since $\eta \ll \sigma$.
Therefore, we have 
        $$e_R((N_R(a) \cap (X \setminus V(M_0)), N_R(b) \cap (X \setminus V(M_0)) \ge \frac{\sigma}{2}|X|^2>0.$$
  Choose an edge $cd\in E_R(N_R(a) \cap (X \setminus V(M_0)), N_R(b) \cap (X \setminus V(M_0)))$. Then we have a 4-cycle $C=(a,c,d,b,a)$ containing exactly three edges in $R$ and $V(C) \cap V(M_0) \subseteq \{a,b\}$. This leads to a contradiction with Claim \ref{even-cycle}.

        If there is a set $A\subseteq X$ with         $(\frac 12-21\sigma)|X| \le |A| \le (\frac 12+21\sigma)|X|$ and $e_R(A,X\setminus A) \le 4\sigma |X|^2\le 4\sigma n^2$ (Lemma~\ref{Class} (2)), then
        $$
       e_R(A,\overline{A})=  e_R(A,X \setminus A) +e_R(A,Y)+e_R(A,Z) \le 4\sigma n^2+\alpha n^2+ \alpha n^2 \le 5 \sigma n^2,
        $$
        where the second inequality holds as $e_{R}(A,Y) \le e_G(X,Y) \le \alpha n^2$, and $e_{R}(A,Z) \le |Z|n \le \alpha n^2$, the last inequality holds since $\alpha \ll \sigma$.
        %{\color{blue}前两个值分别来自两次应用“分类引理”，最后一个是因为$|Z|$足够小。}
        However, according to Lemma~\ref{R-lem} (c),
        $$
        e_R(A,\overline{A}) \ge \frac{1}{2}e_G(A,X \setminus A) \ge\frac{1}{2}(\delta(X)-|A|)|A| \ge (\frac{1}{32}-4\sigma)n^2, 
        $$
where the last inequality holds since $\delta(X) \ge (\frac12-2\eta)n \ge (1-5\eta)|X|$, and $|X|\ge(\frac 12-\alpha)n$.
%, and  $\frac{1}{2}(\delta(X)-|A|)|A|\ge \frac 12\left((1-5\eta)|X|-(\frac 12+21\sigma)|X|\right)(\frac 12+21\sigma)|X|\ge (\frac{1}{32}-4\sigma)n^2$
This leads to a contradiction with $\alpha n^2\ge (\frac{1}{32}-4\sigma)n^2$.

       Finally, suppose there exists a set $A$ satisfying that  $(\frac 12-25\sigma)|X| \le |A| \le (\frac 12+25\sigma)|X|$, 
       $\delta(R[A, X\setminus A])\ge \frac{|X|}5 %\ge \frac n{11} 
       $, $e_R(A,X\setminus A) \ge (\frac{1}{4}-5\sigma)|X|^2 %\ge (\frac{1}{16}-2\sigma)n^2
       $, 
       and $e(R[A]) \le 6\sigma |X|^2 %\le 2\sigma n^2
       $ (Lemma~\ref{Class} (3)).
       %the bipartite graph induced by the edges in $E(R)$ between $A$ and $X \setminus A$ has minimal degree at least ${|X|}/{5} \ge n/11$, 
    %   and
    %    $$e_R(A,X\setminus A) \ge (\frac{1}{4}-5\sigma)|X|^2 \ge (\frac{1}{16}-2\sigma)n^2, e(R[A]) \le 6\sigma |X|^2 \le 2\sigma n^2.$$
        We show that $R[X]=G[A,X\setminus A]$, i.e., $R[X]$ forms a cut of $G[X]$.  %Denote $B=X\setminus A$.
                %Then we can apply Lemma~\ref{bidiameter}~(2) to bipartite graph $R[A,X\setminus A]$ for $\varepsilon=5\sigma$. 
    %    Since $\delta(R[A,X\setminus A]) \ge \frac{|X|}{5} \ge  2\sqrt{5\sigma}|X|+1$ and $e_R(A,X\setminus A) \ge (\frac{1}{4}-5\sigma)|X|^2$.

        First, we claim that $R'=R[A, X\setminus A]$ is $(\sqrt{5\sigma}n,4)$-connected. Choose $U \subseteq V(R')$ with $|U|\le \sqrt{5\sigma}n$. 
    We claim that $\diam(R'-U)\le 4$.
   It is sufficient to show that $\dist(a, b)\le 3$ for any $a \in A\setminus U$ and $b \in (X\setminus A)\setminus U$. Note that $\delta(R')=\delta(R[A, X\setminus A]) \ge \frac{|X|}{5} \ge  2\sqrt{5\sigma}|X|+1$. We have $\min\{|N_{R'}(a)\setminus U|,|N_{R'}(b)\setminus U|\} \ge \delta(R') -|U| \ge \sqrt{5\sigma}n+1$. 
  This implies that $E_{R'}(N_{R'}(a)\setminus U, N_{R'}(b)\setminus U)\neq \emptyset$; otherwise, $e(R')=e_{R'}(A,X\setminus A) \le |A|\cdot|X\setminus A|-|N_{R'}(a)\setminus U|\cdot|N_{R'}(b)\setminus U| \le \frac14 n^2 -(\sqrt{5\sigma} n+1)^2<(\frac14 - 5\sigma) n^2$, leading to a contradiction with $e(R') \ge (\frac14 - 5\sigma) n^2$. 
 Therefore, there exists a path of length at most three connecting $a$ and $b$ in $R'$, we are done.
        %, and hence has a diameter of at most $4$.

        Second, we assert that  $E_R(A,X\setminus A)= E_G(A,X\setminus A)$.
        If not, then there is an edge $ab \in E_G(A,X\setminus A)\setminus E(R)$ with $a \in A$ and $b \in X \setminus A$. 
        Since $R[A,X\setminus A]$ is $(\sqrt{5\sigma}n,4)$-connected and $|M_0| =3$, we have a path $P$ of length at most 4 in $$R[\left(A\setminus V(M_0)\right)\cup \{a,b\},\left(X\setminus(A\cup V(M_0))\right)\cup \{a,b\}]$$ connecting $a$ and $b$. Additionally, $e(P)$ is odd since $a \in A$ and $b\in X\setminus A$. 
      %  {\color{brown}then we can find a path $P$ in it connecting $a$ and $b$ with at most $4$ edges by the diameter of $R[A,X\setminus A]$}, and $e(P)$ is odd by $a \in A$ and $b\in X\setminus A$. 
%  \iffalse
% Given that $|N_R(a) \cap (X\setminus A)|,|N_R(b) \cap A|\ge \frac{|X|}5,$ we must have 
%        $$e_R(N_R(a) \cap (X\setminus A),N_R(b) \cap A) > 0,$$
 %       otherwise, we would deduce that $e_R(A,X\setminus A) \le \frac{|X|^2}4-|N_R(a) \cap (X\setminus A)|\cdot|N_R(b) \cap A| \le \frac{|X|^2}4 -\frac{|X|^2}{25}< (\frac{1}{4}-5\sigma)|X|^2$, a contradiction.
 % By selecting an edge $cd\in E_R(N_R(a) \cap (X\setminus A),N_R(b)\cap A)$ with $c\in N_R(a) \cap (X\setminus A)$ and $d \in N_R(b) \cap A$,
 % \fi
Therefore,   we obtain an even cycle $P+ab$ with an odd number of edges in $E(R)$, and its intersection with $V(M_0)$ is included in $\{a,b\}$,
which contradicts Claim \ref{even-cycle}. 
    %    Therefore, we have $E_R(A,X\setminus A) = E_G(A,X\setminus A)$.
    
Next, we assert that $E(R[A]) \cup E(R[X\setminus A]) = \emptyset$. If not, choose an edge $cd \in E(R[A])\cup E(R[X\setminus A])$. 
     Clearly, 
        $$
       e(G[A])=
       \frac{1}{2}\sum_{u \in A}|N_G(u) \cap A| \ge
       \frac{1}{2}|A|\left(\left(\frac{1}{2}-2\eta\right)n-|X\setminus A|\right) \ge \left(\frac{1}{32}-4\sigma\right)n^2.
        $$ 
Since $e(R[A]) \le 6\sigma |X|^2 \le 2\sigma n^2$, we have $$e(G[A])-e(R[A]) \ge \left(\frac{1}{32}-6\sigma\right)n^2>5|A|>e_G\left(\{c,d\}\cup (V(M_0)\cap X), A\right).$$ 
Therefore, we can choose an edge $uv \in E(G[A])\setminus E(R[A])$ such that $\{u,v\}\cap(\{c,d\}\cup V(M_0))=\emptyset$.
Recall that the bipartite subgraph $R[A,X\setminus A]$ is $(\sqrt{5\sigma}n,4)$-connected. According to Lemma~\ref{bidiameter}, for the two pairs $(u,c)$ and $(v,d)$, we can construct two disjoint paths, whose inner vertices avoid $V(M_0)$, $P_{uc}$ (connecting $u$ and $c$) and $P_{vd}$ (connecting $v$ and $d$) such that $e(P_{uc})+e(P_{vd})\le 8$, and $e(P_{uc})$ and $e(P_{vd})$ have the same parity since $\{u,v\} \subseteq A$ and $\{c,d\} \subseteq A \text{ or } X\setminus A$. Therefore, we have an even cycle $P_{uc}+cd+P_{vd}+vu$ with  an odd number of edges in $E(R)$ and at most two vertices in $V(M_0)$ (should belong to $\{c,d\}$),  a contrary  to Claim \ref{even-cycle}.

In summary, we have $R[X]=G[A,X\setminus A]$.

\noindent \textbf{Step II.} We prove that there exists a partition $Z=Z_A \cup Z_{X\setminus A}$ such that $R[X \cup Z]=G[A\cup Z_A,(X\setminus A)\cup Z_{X\setminus A}]$.
       
         Consider a random vertex $z\in Z$. We claim that either $N_R(z,A)=N_G(z,A)$ and $ N_R(z,X\setminus A)=\emptyset$, or $N_R(z,A)=\emptyset$
        and $N_R(z, X\setminus A)=N_G(z,X\setminus A)$. Recall that $\delta_X(Z) \ge \frac34 \eta n$, without loss of generality, we may assume  $\delta_A(Z) \ge \frac38 \eta n$ by pigeonhole principle. Then $|N_G(z) \cap (A\setminus V(M_0))|\ge \frac38 \eta n-3>0$.
 
       % If there exist $u_1 \in N_R(z) \cap A$ and $u_2 \in (N_{G \setminus R}(z) \cap A)\setminus V(M_0)$, 
     Fix a vertex $u_1 \in N_G(z) \cap (A\setminus V(M_0))$.  For any $u_2\in N_G(z) \cap X$,
        choose another vertex $u_3 \in (N_G(u_2) \cap (X\setminus A))\setminus V(M_0)$.
       Recall that $R[A,X\setminus A]$ is $(\sqrt{5\sigma}n,4)$-connected and $|V(M_0)| =6 \le \sqrt{5\sigma}n$. We have  $R[A,X\setminus A]-V(M_0)$ is $(\sqrt{5\sigma}n-6,4)$-connected.
       Then we can find a path $P_{u_1u_3}$ in $R[A,X\setminus A]-V(M_0)$ connecting $u_1$ and $u_3$ with length at most $4$.
      Since $u_1 \in A$ and $u_3 \in X\setminus A$, we have $e(P_{u_1u_3})$ is odd. Hence $C=zu_2u_3+P_{u_1u_3}+u_1z$ is an even cycle in $G[X\cup Z]$ of length at most $6$, and, additionally, $C$ contains at most one vertex (possibly $u_2$) in $M_0$.
      According to Claim~\ref{even-cycle}, $|\{zu_1,zu_2,u_2u_3\} \cap E(R)|$ must be odd, implying that $|\{zu_1,zu_2,u_2u_3\} \cap E(R)|=1$ or 3. 
      If $u_2 \in A$, then by Step I, $u_2u_3\in E(R)$. Consequently,  we have either $\{zu_1,zu_2\} \subseteq E(R)$ or $\{zu_1,zu_2\} \cap E(R) = \emptyset$. 
    If $u_2 \subseteq X\setminus A$, then by Step I,  $u_2u_3\notin E(R)$. In this case, either $zu_1 \in E(R)$ and $zu_2 \notin E(R)$ or vice versa.
Therefore, when $u_2$ traverses all possible selections within $N_G(z) \cap X$, we deduce that either $E_G(z,A)=E_R(z,A)$ and $E_R(z,X\setminus A)=\emptyset$, or $E_R(z,A)=\emptyset$ and $E_G(z,X\setminus A)=E_R(z,X\setminus A)$.

    Now we can partition $Z$ into
        $$
        Z_A=\{z\in Z: E_R(z,A)=\emptyset \text{ and } E_G(z,X\setminus A)=E_R(z,X\setminus A)\},
        $$
        and
        $$
        Z_{X\setminus A}=\{z\in Z: E_R(z,X\setminus A)=\emptyset \text{ and } E_G(z,A)=E_R(z,A)\}.
        $$
    To show $R[X\cup Z]$ forms a cut of $G[X\cup Z]$, it is sufficient to prove that $E(R[Z_A])\cup E(R[Z_{X\setminus A}])=\emptyset$ and $R[Z_A,Z_{X\setminus A}]=G[Z_A,Z_{X\setminus A}]$. 
    For any edge $z_1z_2 \in E(G)$  with $z_1\in Z$ and $z_2 \in Z\setminus V(M_0)$. Choose $u_1 \in (N_G(z_1) \cap A)\setminus V(M_0)$ and $u_2 \in (N_G(z_2)\cap (X\setminus A))\setminus V(M_0)$. By the definition of $Z_A$ and $Z_{X\setminus A}$, $|\{z_1u_1,z_2u_2\}\cap E(R)|=1$ if $\{z_1,z_2\} \in \binom{Z_A}{2}\cup \binom{Z_{X\setminus A}}{2}$, and $|\{z_1u_1,z_2u_2\}\cap E(R)|=0$ or 2, otherwise.
        Since $R[A,X\setminus A]-V(M_0)$ is $(\sqrt{5\sigma}n-6,4)$-connected, we can find a path $P_{u_1u_2}$ connecting $u_1$ and $u_2$ of length at most 4 in it. Additionally, $|E(P_{u_1u_2})|$ is odd since $u_1\in A$ and $u_2 \in X\setminus A$. Therefore, we have an even cycle $C=z_1u_1P_{u_1u_2}u_2z_2z_1$ containing at most one vertex (possible $z_1$) in $V(M_0)$.  
        By Claim~\ref{even-cycle}, $|\{z_1z_2,z_1u_1,z_2u_2\}\cap E(R)|$ must be odd. 
        Therefore, if $\{z_1,z_2\} \in \binom{Z_A}{2}\cup \binom{Z_{X\setminus A}}{2}$, we have $z_1z_2 \notin E(R)$, and $z_1z_2\in E(R)$, otherwise.
 This  deduces that $E(R[Z_A])\cup E(R[Z_{X\setminus A}])=\emptyset$ and $E(R(Z_A,Z_{X\setminus A}))=E(G(Z_A,Z_{X\setminus A}))$. We are done.

 %       Overall, we have proved that $E_G(A\cup Z_A) \cup E_G((X\setminus A)\cup Z_{X\setminus A}) \subseteq E(G)\setminus E(R)$ and $E_G(A\cup Z_A,(X\setminus A)\cup Z_{X\setminus A}) \subseteq E(R)$,
%        this implies that $R[X\cup Z]$ forms a cut of $G[X\cup Z]$, and $A \cup Z_A$, $(X\setminus A) \cup Z_{X\setminus A}$ are two parts of $R[X\cup Z]$.
    \end{proof}
    
    \begin{claim} \label{2-cut}
    If $|M_0|\ge 3$, then it is not possible for both $R[X]$ to form a cut of $G[X]$ and $R[Y]$ to form a cut  of $G[Y]$, simultaneously.
    
%     If $|M_0|\ge 3$ and {\color{red} $R[X]$ and $R[Y]$ both forms a cut.} Then $R$ forms a cut of $G$.
    
   %     If $M_0$ is a matching between $X \cup Z$ and $Y$, $|E(M_0)|\ge 3$, $R[X],R[Y]$ both forms a cut. Then $R$ forms a cut of $G$.
    \end{claim}

    \begin{proof}
Suppose this is the case. According to Claim~\ref{claim2}, there exist $A \subseteq X$ and $C\subseteq Y$ with $(\frac 12-25\sigma)|X| \le |A| \le (\frac 12+25\sigma)|X|$ and $(\frac 12-25\sigma)|Y| \le |C| \le (\frac 12+25\sigma)|Y|$, and a partition $Z=Z_A \cup Z_{X\setminus A}$ such that $R[X \cup Z]=G[A\cup Z_A,(X\setminus A)\cup Z_{X\setminus A}]$ and $R[Y]=G[C,Y\setminus C]$.  We denote $X\setminus A$ as $B$, and $Y\setminus C$ as $D$, respectively.

%denote $R[X]$ as $E(A,B)$ for the partition $X=A\cup B$ and $R[Y]$ as $E(C, D)$ for the partition $Y=C\cup D$.
     %Let $A,B$ be two parts of $R[X]$ and $C,D$ be two parts of $R[Y]$.
     %   By condition we have $E_R(A)=E_R(B)=E_R(C)=E_R(D)=\emptyset$

        Randomly choose two non-adjacent edges $x_1y_1, x_2y_2\in E_G(X \cup Z, Y)$.
      Define a mapping $\phi : (X\cup Z)^2 \cup Y^2\rightarrow\{0,1\}$ such that $\phi(u,v)=0$ if $(u,v)\in (A \cup Z_A)^2\cup(B \cup Z_B)^2 \cup C^2\cup D^2$; otherwise, $\phi(u,v)=1$. Here $S^2$ denotes the Cartesian product of the set $S$. 
      By the definition of $\phi$ and the properties of $A, B, C, D, Z_A$, and $Z_B$,  we  deduce that
      $$\phi(u,v)=\phi(v,u) \text{ for any } (u,v)\in (X\cup Z)^2 \cup Y^2$$  and 
      \begin{align}\label{mod}
          \phi(u,v) = |\{uv\} \cap E(R)| \text{ for any $uv \in E(G)$}.
      \end{align}
      We claim that if $v_1, v_2, \dots, v_{\ell}$ are $\ell$ vertices contained either in $X \cup Z$ or in $Y$, then
\begin{align}\label{sum}
    \sum_{i=1}^{\ell-1} \phi(v_i, v_{i+1}) \equiv \phi(v_1, v_{\ell}) \pmod{2}.
\end{align}
We prove the equation (\ref{sum}) by induction on $\ell$. When $\ell=2$, (\ref{sum}) holds trivially. Now assume $\ell\ge 3$. Without loss of generality, assume $\{v_1, v_2, \dots, v_{\ell}\}\subseteq X\cup Z$ and $v_1\in A\cup Z_A$.  
By the induction hypothesis,  
\begin{align}\label{sum1}
\sum_{i=1}^{\ell-1} \phi(v_i, v_{i+1}) \equiv \sum_{i=1}^{\ell-2} \phi(v_i, v_{i+1})+\phi(v_{\ell-1}, v_{\ell})\equiv \phi(v_1, v_{\ell-1})+ \phi(v_{\ell-1}, v_{\ell})\pmod{2}. 
\end{align}
If  $v_{\ell-1}$ and $v_\ell$ are both in $A\cup Z_A$, then $\phi(v_1, v_{\ell-1})+ \phi(v_{\ell-1}, v_{\ell})= 0+0= 0=\phi(v_1, v_\ell)$. If  $v_{\ell-1}$ and $v_\ell$ are both in $B\cup Z_B$, then $\phi(v_1, v_{\ell-1})+ \phi(v_{\ell-1}, v_{\ell})= 1+0= 1=\phi(v_1, v_\ell)$. If $v_{\ell-1}$ and $v_\ell$ are in different sets, without loss of generality, assume $v_{\ell-1}\in A\cup Z_A$ and $v_\ell\in B\cup Z_B$, then $\phi(v_1, v_{\ell-1})+ \phi(v_{\ell-1}, v_{\ell})= 0+1= 1=\phi(v_1, v_\ell)$. 
Combining with (\ref{sum1}), we obtain (\ref{sum}).
%suffices to prove this for the case $\ell = 3$, which can be checked by considering all possible configurations of $v_1, v_2, v_3$. This can be done by applying the result to $v_1$ and each pair of consecutive vertices (other than $v_1$) in the sequence, which then leads to the conclusion in equation (\ref{sum}).

        Next choose vertices
        $a \in A \cap N_G(x_1)$, $b\in B\cap N_G(x_2)$, $c\in C\cap N_G(y_1)$ and $d\in D\cap N_G(y_2)$, and edges $a'b' \in E_R(A\setminus \{a,x_1,x_2\}, B\setminus \{b,x_1,x_2\})$ and $c'd' \in E_R(C\setminus \{c,y_1,y_2\}, D\setminus \{d,y_1,y_2\})$
    (the existence of edges $a'b'$ and $c'd'$ is guaranteed by $\min\{\delta(Y), \delta(X) \}\ge \left(\frac{1}{2}- 2\eta\right)n$).

  Recall that $x_i^1,x_i^2$ are two chosen neighbors of $z_i \in Z$ in $X$. Since $\delta_X(Z)\ge \frac 34\eta n$ is large enough, without loss of generality, we may assume that $x_i^1,x_i^2 \in A\setminus \{x_1,x_2,a,a'\}$ for all $z_i \in Z\setminus \{x_1,x_2\}$.
   Since
     \begin{eqnarray*}
         \delta(G[A\setminus \{x_1,x_2\}]) &\ge& \delta(X)-|X\setminus A|-2\\
                                           &=&|A|-(|X|-\delta(X)+2)\nonumber\\ &\ge& |A|-5\eta|X|\nonumber-2 \\ &\ge& |A|-\frac{5\eta}{\frac{1}{2}-25\sigma}|A|-2\nonumber \\ 
                                           &\ge& \left(\frac{1}{2}-11\eta\right)|A\setminus\{x_1,x_2\}|,
     \end{eqnarray*} we may apply Lemma \ref{Fixpair} (1) to $G[A\setminus \{x_1,x_2\}]$ with parameters $\varepsilon_1=11\eta$, $\varepsilon_2=2\alpha$, and $\ell$  pairs of vertices (where $\ell \le |Z|+1  \le 2\alpha n$)
     $$\{(a,a')\}\cup \bigcup_{i:z_i \in Z\setminus \{x_1,x_2\}}\{(x_i^1,x_i^2)\}.$$ 
    This yields a path $P_{aa'}$ from $a$ to $a'$ that covers exactly $(A \cup Z)\setminus \{x_1,x_2\}$ (where the path $P_{aa'}$ can be obtained from the resulting Hamilton cycle containing the edges $aa'$, $x_i^1x_i^2$ for $\{i : z_i \in Z\setminus \{x_1,x_2\}\}$ in  $G[A\setminus\{x_1,x_2\}]+\{aa'\}+\bigcup_{i:z_i \in Z\setminus \{x_1,x_2\}}\{x_i^1x_i^2\}$ by replacing each edge $x_i^1x_i^2$ by $x_i^1z_ix_i^2$ and deleting the edge $aa'$). Similar path constructions $P_{vv'}$ for $v \in \{b,c,d\}$ covering their corresponding vertex parts can be executed without requiring vertex embeddings in $Z$. 
Let $P_{x_1x_2} = x_1a + P_{aa'} + a'b' + P_{bb'} + bx_2$ and $P_{y_1y_2} = y_2d + P_{dd'} + d'c' + P_{cc'} + cy_1$.
Therefore, the cycle
$$
C_0 = P_{x_1x_2} + x_1y_1 + P_{y_1y_2} + x_2y_2
$$
is a Hamiltonian cycle. Note that $P_{x_1x_2}$ is contained in $G[X \cup Z]$ and $P_{y_1y_2}$ is contained in $G[Y]$.
Applying equation (\ref{mod}) to all edges in $E(P_{x_1x_2})$ and using equation (\ref{sum}), we have
\begin{align*}
    |E(P_{x_1x_2}) \cap E(R)| 
    &\equiv \sum_{e = uv \in E(P_{x_1x_2})} \phi(u, v) \pmod{2} \\
    &\equiv \phi(x_1, x_2) \pmod{2}.
\end{align*}
Similarly, for $P_{y_1y_2}$, we obtain
\begin{align*}
    |E(P_{y_1y_2}) \cap E(R)| 
    &\equiv \phi(y_1, y_2) \pmod{2}.
\end{align*}
Recall that $C_0$ contains an even number of edges in $E(R)$. By combining the two equations above, we conclude that
\begin{align}\label{2edges}
    |E(R) \cap \{x_1y_1, x_2y_2\}| 
    \equiv \phi(x_1, x_2) + \phi(y_1, y_2) \pmod{2}.
\end{align}

Fix  three  edges $e_1=u_1w_1$, $e_2=u_2w_2$, and $e_3=u_3w_3$ in $M_0$.
%in $E_G(X \cup Z, Y)$ since $|M_0| \ge 3$ ($u_i \in X \cup Z$ and $v_i \in Y$).
Without loss of generality, we assume that $u_1 \in A$, $w_1 \in C$ and assume $u_1w_1 \in E(R)$. For $i \in \{1, 2, 3\}$, let
\begin{align*}
    E_i &= \{ e = uv \in E_G(X \cup Z, Y) : \{u, v\} \cap \{u_i, w_i\} = \emptyset \}.
\end{align*}
We assert that $E_G(X \cup Z, Y) = E_1 \cup E_2 \cup E_3$. Otherwise, there exists an edge $e_0$ not in any of these three sets. Then $e_0$ intersects all three edges $e_1, e_2$ and $e_3$, which is impossible.

For any edge $uw \in E_G(X \cup Z, Y)$, we claim that
\begin{align} \label{E_G(Xcup Z, Y)}
    |E(R) \cap \{u_1w_1, uw\}| 
    &\equiv \phi(u_1, u) + \phi(w_1, w) \pmod{2}.
\end{align}
If $uv \in E_1$, then $\{u, w\} \cap \{u_1, w_1\} = \emptyset$. According to (\ref{2edges}), we have (\ref{E_G(Xcup Z, Y)}).
 In particular, for $u_2w_2\in E_1$, we have 
\begin{align} \label{E_1E_2}
    |E(R) \cap \{u_1w_1, u_2w_2\}| 
    &\equiv \phi(u_1, u_2) + \phi(w_1, w_2) \pmod{2}.
\end{align}
If $uv \notin E_1$, then either $uv \in E_2$ or $uv \in E_3$. Without loss of generality, assume that $uv \in E_2$. This implies that $\{u, w\} \cap \{u_2, w_2\} = \emptyset$. According to (\ref{2edges}), we have 
\begin{align*}
    |E(R) \cap \{u_2w_2, uw\}| 
    &\equiv \phi(u_2, u) + \phi(w_2, w) \pmod{2}.
\end{align*}
By summing  the corresponding left-   and right-hand sides of this equation and equation (\ref{E_1E_2}), we obtain (\ref{E_G(Xcup Z, Y)}).

If $uw \in E_G(A\cup Z_A,C)$, then by (\ref{E_G(Xcup Z, Y)}), $|E(R) \cap \{u_1w_1, uw\}| \equiv \phi(u_1, u) + \phi(w_1, w)=0+0\equiv 0\pmod2$. This implies that $uv\in E(R)$ as $u_1w_1 \in E(R)$. Similarly, we have $uv\in E(R)$ if $uv \in E_G(B\cup Z_B,D)$, and $uv \in E(G)\setminus E(R)$ if $uv \in E_G(B\cup Z_B,C) \cup E_G(A\cup Z_A,D)$. It follows that $E_R(A\cup Z_A, D) \cup E_R(B\cup Z_B, C) = \emptyset$ and $E_R(A\cup Z_A, C) \cup E_R(B\cup Z_B,D) = E_G(A\cup Z_A ,C) \cup E_G(B\cup Z_B, D)$. 
Combining with the fact that $R[X \cup Z]=G[A\cup Z_A, B\cup Z_B]$ and $R[Y]=G[C, D]$, 
we finally have $E(R)=E_G(A\cup Z_A\cup D, B\cup Z_B\cup C)$, i.e., $R$ forms a cut of $G$. This leads to a contradiction with the property (c) of $R$. 

%By the definition of $A,B,C,D$ we finally have 
%$$
%E_R(A \cup D) \cup E_R(B \cup C) = \emptyset \text{ and } E_R(A \cup D,B \cup C)= E_G(A \cup D,B \cup C),
%$$
%i.e. $R$ forms a cut with $A\cup D$ and $B\cup C$ as two parts.

        %{\color{red} The number of edges in $(E(C_0) \cap E(R)) \setminus \{x_1y_1,x_2y_2\}$ satisfies $|(E(C_0) \cap E(R)) \setminus \{x_1y_1,x_2y_2\}|\equiv \phi(x_1,a)+\phi(x_2,b)+\phi(y_1,c)+\phi(y_2,d)\pmod 2$. Since $C_0$ must contain an even number of edges in $E(R)$, $|E(R) \cap \{x_1y_1,x_2y_2\}| \equiv \phi(x_1, x_2) + \phi(y_1,y_2)\pmod 2$. And for any $e \in E_G(X \cup Z,Y)$ there is another edge $e' \in E_G(X \cup Z,Y)$ disjoint from $e$ as $|E(M_0)| \ge 3$, implying that $R$ is a cut of $G$. \textbf{NEED A REASON?}}
    \end{proof}

    \begin{claim} \label{0-cut}
    If $|M_0| \ge 3$, then it is not possible for both $R[X]=G[X]$ and $R[Y]=G[Y]$, simultaneously.
    %{\color{red}$M_0$ is a matching between $X,Y$, $|E(M_0)|\ge 3$,??? then $R=G$.}
    \end{claim}

    \begin{proof}
    Suppose this is the case.
     Choose three edges in $M_0$, denoted as $u_iy_i$ with $u_i\in X\cup Z$ and $y_i\in Y \text{ for } i\in [3]$. Choose $x_i \in N_G(u_i) \cap X$ for $i \in [3]$, and let $X'=X\setminus \{u_1,u_2,u_3\}$ and $Z'=Z\setminus \{u_1,u_2,u_3\}$. 
    Again since $\delta_X(Z)\ge \frac 34\eta n$ is large enough, we may assume that $x_i^1,x_i^2 \in X'\setminus \{x_1,x_2,x_3\}$ for all $1 \le i \le |Z|$.
Since  $\delta(G[X'])\ge (\frac12-2\eta)n-2\ge (1-5\eta)|X|-2 \ge (1-6\eta)|X'|$, we may apply Lemma \ref{Fixpair} (1) to $G[X']$ with parameters $\varepsilon_1=6\eta$, $\varepsilon_2=2\alpha$ and $|Z'|+1$ pairs of vertices (where  $|Z'|+1 \le 2\alpha n$)
        $$
        \{(x_j,x_k)\}\cup \left(\bigcup_{i: z_i \in Z'} \{(x_i^1x_i^2)\}\right) \text{  for $\{j,k\} \in \binom{[3]}{2}$}.
        $$
This yields an $x_j$-$x_k$- path $P_{x_jx_k}$  that covers exactly $(X \cup Z)\setminus \{u_j,u_k\}$ (where $P_{x_jx_k}$ can be obtained from the resulting Hamilton cycle  containing edges $x_jx_k$ and $x_i^1x_i^2$ for $i\in\{i : z_i\in Z'\}$ in  $G[X']+\{x_1x_2\}+\bigcup_{i : z_i\in Z'}\{x_i^1x_i^2\}$ by replacing each $x_i^1x_i^2$ by $x_i^1z_ix_i^2$ and deleting the edge $x_1x_2$). By the choice of $x_i^1,x_i^2$, we have that $E(P_{x_jx_k})$ contains an even number of edges in $E(G)\setminus E(R)$. 
%So we have constructed three paths $P_{x_1x_2},P_{x_2x_3},P_{x_3x_1}$. 
Similar path constructions $P_{y_jy_k}$ for $\{j,k\} \in \binom{[3]}{2}$ covering exactly $Y$ can be executed without requiring the embedding of vertices in $Z$.
Therefore, for any $\{j,k\} \in \binom{[3]}{2}$, the cycle
        $$
        C_{jk}=x_ju_jy_j+P_{y_jy_k}+y_ku_kx_k+P_{x_j,x_k}
        $$
        is a Hamiltonian cycle of $G$, 
        where $E(P_{y_jy_k}) \subseteq E(R)$ and the path $E(P_{x_jx_k})$ contains an even number of edges in $E(G)\setminus E(R)$. 
         %   This implies that $|(E(P_{y_jy_k})\cup E(P_{x_jx_k})\setminus E(R)|$ is even. 
        By the construction of $R$,    we know that $|E(C_{jk}) \setminus E(R)|$ is odd.  
This implies that
$|\{u_jy_j,u_jx_j,u_ky_k,u_kx_k\} \cap E(R)|$ must be odd. Let $|\{u_ky_k,u_kx_k\} \cap E(R)|=h_k$ and $|\{u_jy_j,u_jx_j\} \cap E(R)|=h_j$, then 
$$|\{u_jy_j,u_jx_j,u_ky_k,u_kx_k\} \cap E(R)|=h_k+h_j.$$
However, the sum of it over all $\{j,k\}\in \binom{[3]}{2}$
$$\sum_{\{j,k\}\in \binom{[3]}{2}}(h_j+h_k)\equiv 0 \pmod 2,$$ which is impossible.

   %{\color{red}     which is impossible since the sum of all $\{j,k\}\in \binom{[3]}{2}$, which is odd in this case, should be even.}
       
  %      {\color{brown}Now for any $\{j,k\} \in \binom{[3]}{2}$, we can construct the following Hamilton cycle:
   %     $$
    %    C_{jk}=x_ju_jy_jP_{y_j,y_k}y_ku_kx_kP_{x_j,x_k}x_j,
     %   $$
      %  where $E(P_{y_jy_k}) \subseteq E(R)$ and $E(P_{x_jx_k})$ contains an even number of edges in $E(G)\setminus E(R)$, this implies that $|%(E(P_{y_jy_k})\cup E(P_{x_jx_k})\setminus E(R)|$ is even. And we know that $|E(C_{jk}) \setminus E(R)|$ is odd by the definition of $R$,}
    %    so we can deduce that $|\{u_jy_j,u_jx_j,u_ky_k,u_kx_k\} \cap E(R)|$ is odd,
     %   which is impossible since the sum of all $\{j,k\}\in \binom{[3]}{2}$, which is odd in this case, should be even.

    \end{proof}

  Since the inequality
        \begin{eqnarray}\label{3match}
            |Z|+m(X,Y) \ge \frac{3}{4}(\frac{4}{3}|Z|+m(X,Y)) \ge \frac52
        \end{eqnarray}
     holds       and recalling that $\delta_X(Z), \delta_Y(Z) \ge \frac34 \eta n$,
        we conclude that there exists is a three-edge matching between $X$ and $Y\cup Z$, and as well as one between $X\cup Z$ and $Y$. 
According to Claim \ref{claim2},
        either $R[X]=G[X]$ or $R[X]$ forms a cut of $G[X]$, and $R[Y]=G[Y]$ or $R[Y]$ forms a cut of $G[Y]$.
        Claim \ref{2-cut} addresses the scenario where both
        $R[X],R[Y]$ are cuts, while Claim \ref{0-cut} considers the case where $R[X]=G[X]$ and $R[Y]=G[Y]$.
        Now we only need to discuss the case of $G[X]=R[X]$ and $R[Y]$ is a cut of $G[Y]$.
        
        \begin{claim}
            If $\frac{4}{3}|Z|+m(X,Y) \ge \frac{10}3$, then it is not possible for both $R[X]=G[X]$ and $R[Y]$ is a cut of $G[Y]$, simultaneously.
        \end{claim}
        \begin{proof}
            Suppose this is the  case. Then by Claim~\ref{claim2}, there exists a subset $C$ of $Y$ such that
            $(\frac 12-25\sigma)|Y| \le |C| \le (\frac 12+25\sigma)|Y|$ and $R[C,Y\setminus C]$ is $(\sqrt{5\sigma}n,4)$-connected. We assert that at least one of the following holds:
            \begin{itemize}
                \item[(A1)] There exist distinct vertices: $w \in Z$ and $x_1,x_2 \in X \cup Z$ and $y_1,y_2\in Y$ such that $x_1y_1,x_2y_2\in E(G)$;
                \item[(A2)] There exist $4$ disjoint edges $e_1,e_2,e_3,e_4$ in $G[X,Y]$.
            \end{itemize}
            
            If $|Z| \ge 1$, we choose some $w \in Z$, then by (\ref{3match}) we know there exist two disjoint edges in $G[(X\cup Z)\setminus \{w\},Y]$, denoted as $x_1y_1$ and $x_2y_2$ ($x_i \in X\cup Z$ and $y_i \in Y$).
            %And by $\delta_X(Z),\delta_Y(Z) \ge \frac34 \eta n$ we can find four neighbors of $w$ labelled as $\{x_3,x_4,y_3,y_4\}$, such that $\{x_3,x_4\}$ are in $X\setminus \{x_1,x_2\}$ and $\{y_3,y_4\}$ are in $Y\setminus \{y_1,y_2\}$. 
            This corresponds to (A1). If $|Z|=0$, then by (\ref{cod:biclique}), we have $m(X,Y) \ge 4$, this is the case of (A2).

            Now, we choose $x_3,x_4 \in (N_G(w) \cap X)\setminus \{x_1,x_2\}$ and $y_3, y_4 \in (N_G(w) \cap Y)\setminus \{y_1,y_2\}$ in the case (A1). In the case (A2), label $e_i$ as $x_iy_i$ for $i \in [4]$, where $x_i \in X$ and $y_i \in Y$. 
            %Note these notations won't cause ambiguity.
             No matter in the case of (A1) or (A2), we fix $x_i$'s and $y_i$'s for $i \in [4]$. Without loss of generality, we set $y_3 \in C$.
            We are going to construct some paths between $x_1,x_3$ and $y_1,y_3$, respectively. 
            Choose 
            %$x_3' \in (N_G(x_3) \cap X)\setminus \{x_1,x_2,x_3,x_4\}$ and 
            $y_3' \in (N_G(y_3) \cap C)\setminus \{y_1,y_2,y_3,y_4\}$ (hence edge $y_3y_3'$ is in $E(G)\setminus E(R)$ by $\{y_3,y_3'\}\subseteq C$). 
            
            Let $X_1=(X \cup Z)\setminus \{x_2,x_4\}$.
            Recall that Claim~\ref{connected}
            tells us that $G[X]$ is $(\frac12 \eta n,2)$-connected  and $G[X\cup Z]$ is $(\frac12 \eta n,4)$-connected. Then $G[X_1\setminus Z]$ is $(\frac12 \eta n-2,2)$-connected and $G[X_1]$ is $(\frac12 \eta n-2,4)$-connected. 
        Applying Claim~\ref{CL: differlengthpath} to $G[X_1\setminus Z]$ and $G[X_1]$, we obtain two $x_1$-$x_3$-paths $P_{13}$ and $P_{13}'$ such that $|E(P_{13})|=4$ and $|E(P_{13}')|=3$ in $G[X_1]$.
        
               On the other hand, let $Y_1=Y\setminus \{y_2,y_3,y_4\}$.
            Recall that $R[Y]=R[C,Y\setminus C]$ is $(\sqrt{5\sigma}n,4)$-connected. We have $R[Y_1]$ is also bipartite and $(\sqrt{5\sigma}n-3,4)$-connected. Hence,
            we can find two paths of length at most $4$: a $y_3$-$y_1$-path $P_{y_3y_1}$ in $R[Y_1]$ and a $y'_3$-$y_1$ path $P_{y_3'y_1}$ in $R[Y_1\setminus \{y_3\}]$ (which may intersect).
           Consequently, it is evident both sets of edges from $P_{y_3y_1}$ and $P_{y_3'y_1}$ are contained in $E(R)$, and they share the same parity.
           Without loss of generality, assume $|E(P_{y_3y_1})|$ is even.
           Define $Q_{13}=P_{y_3y_1}$ and $Q_{13}'=P_{y_3'y_1}+y'_3y_3$. Then $|E_G(Q_{13})|$ be even and $|E_G(Q_{13}')|$ be odd. However, $|E(Q_{13}) \cap E(R)|$ and $|E(Q_{13}')\cap E(R)|$ have the same parity.

            Overall, we have identified paths \( P_{13}, P_{13}', Q_{13}, Q_{13}' \) such that:  
\begin{itemize}
    \item \( |E_G(P_{13})| = |E_R(P_{13})| \) is even;
    \item \( |E_G(P_{13}')| = |E_R(P_{13}')| \) is odd;
    \item \( |E_G(Q_{13})| \) is even, and \( |E_G(Q_{13}')| \) is odd;
    \item \( |E_G(Q_{13}) \cap E(R)| \) and \( |E_G(Q_{13}') \cap E(R)| \) are either both even or both odd;
    \item all of these paths are disjoint from \( \{x_2, y_2, x_4, y_4\} \) and have length at most 4.
\end{itemize}

We will find an $R$-parity-switcher $W$ and an edge in $G[X\cup Z,Y]$ avoiding $V(W)$, in cases (A1) and (A2), respectively.

            If case (A1) holds, consider the cycles 
\[
C_1 =wy_3+Q_{13}'+y_1x_1+P_{13}+x_3w  \text{ and }  C_2 = wy_3+Q_{13}+y_1x_1+P_{13}'+x_3w.
\]
From the construction, it follows that \( |E_G(C_1)| \) and \( |E_G(C_2)| \) are both even. However, the parities of \( |E(C_1) \cap E(R)| \) and \( |E(C_2) \cap E(R)| \) are differ. Therefore, at least one of them must be odd. Without loss of generality, assume \( |E(C_1) \cap E(R)| \) is odd (otherwise, the roles of  $X \cup Z$  and  $Y$  can be swapped in the following argument).

Label \( P_{13} \) as 
$(x_1 =) u_1 u_2 \ldots u_t \ldots u_{2t-1}(=x_3)$,
and \( Q_{13}' \) as 
$(y_3=)v_1 v_2 \ldots v_\ell \ldots v_{2\ell}(=y_1)$,
where \( t \in \{2, 3\} \) and \( \ell \in \{1, 2, 3\} \). Hence, the total number of edges in \( C_1 \) is \( |E(C_1)| = 2t + 2\ell \).
            Let $X_2=(X\cup Z)\setminus\{w,u_t,x_2\}$ and $Y_2=Y\setminus \{y_2,v_\ell\}$. So $G[X_2]$ is $(\frac{1}{2}\eta n-3,4)$-connected and $G[Y_2]$ is $(\sqrt{5\sigma}n-2,4)$-connected.
            %(In sure, we even have $R[Y_2]$ is $(2\sqrt{\sigma}n,4)$-connected, but we won't need to require edges being in $E(R)$ or not). 
            Applying  Lemma~\ref{bidiameter} to $G[X_2]$ for $t-1 \le 2$ pairs of vertices: $(u_1,u_{2t-1}),(u_2,u_{2t-2}),\dots,(u_i,u_{2t-i}),\dots,(u_{t-1},u_{t+1})$,  we obtain $t-1$ disjoint paths $P_1,\dots,P_i,\dots,P_{t-1}$ in $G[X_2]$, where $P_{i}$ connects $u_i$ and $u_{2t-i}$ and has at most $4$ edges. Applying Lemma~\ref{bidiameter} to $G[Y_2]$ for $\ell \le 3$ pairs of vertices: $(y_4,y_1),(v_1,v_{2\ell-1}),(v_2,v_{2\ell-2}),\dots,(v_{\ell-1},v_{\ell+1})$, we obtain $\ell$ disjoint paths $Q_{1},\dots,Q_{i},\dots,Q_{{\ell-1}}$ and $Q_{y_4y_1}$ in $G[Y_2]$, where $Q_{i}$ connects $v_i$ and $v_{2\ell-i}$ and $Q_{y_4y_1}$ connects $y_1$ and $y_4$ and each has at most $4$ edges. Let $Q_{w}=Q_{y_4y_1}+wy_4$. Then $|E(Q_{w})|\le 5$. Set $P_t=u_t$ and $Q_\ell=v_\ell$.
        Therefore, we construct an $R$-parity-switcher 
               $$W=W(R,C_1, P_t, P_{1},\dots,\dots,P_{{t-1}},Q_{w},Q_{1},\dots,\dots,Q_{{\ell-1}}, Q_{\ell}, u_t,v_\ell).$$
        Moreover,  $$|V(W)| \le 2+5(\ell+t-2)+6 \le 28 \le \frac{5}{12}\eta n +100.$$
        This, combined with the fact that the edge $x_2y_2$ is disjoint from $V(W)$, leads to a contradiction to Claim~\ref{noR-parity}.

            If case (A2) holds,  consider the cycles:  
            $$
             C_1 = x_3y_3+Q_{13}+y_1x_1+P_{13}  \text{ and }  C_2 = x_3y_3+Q_{13}'+y_1x_1+P_{13}'.
            $$
            Then both $|E_G(C_1)|$ and $|E_G(C_2)|$ are even. However, the parities of \( |E(C_1) \cap E(R)| \) and \( |E(C_2) \cap E(R)| \) are differ. Thus one of these cycles must have an odd number of edges in \(R\).

        When $|E(C_1) \cap E(R)|$ is odd, we label $P_{13}$ as $(x_1=)u_1u_2\dots u_t \dots u_{2t-1}(=x_3)$ and $Q_{13}$ as $(y_3=)v_1v_2\dots v_\ell \dots v_{2\ell-1}(=y_1)$ for some $t,\ell \in \{2,3\}$. It follows that $|E(C_1)|=2t+2\ell-2$.
Similar to the (A1) case,  we obtain $t-1$ disjoint paths $P_1,\dots,P_i,\dots,P_{t-1}$ in $G[X\setminus\{u_t, x_2\}]$, where $P_{i}$ connects $u_i$ and $u_{2t-i}$ and contains at most $4$ edges. Similarly, there are $\ell-1$ disjoint paths $Q_{1},\dots,Q_{i},\dots,Q_{{\ell-1}}$ in $G[Y\setminus\{v_\ell, y_2\}]$, where $Q_{i}$ connects $v_i$ and $v_{2\ell-i}$ and each has at most $4$ edges. Set $P_t=u_t$ and $Q_\ell=v_\ell$.
%Then we can obtain $t+\ell-2$ disjoint paths of length at most $4$: $P_{1,2t-1},\dots,P_{i,{2t-i}},\dots,P_{{t-1},{t+1}}$ and $Q_{1,2\ell-1},\dots,Q_{i,{2\ell-i}},\dots,Q_{{\ell-1},{\ell+1}}$, where $P_{i,2t-i}$ connects $u_i$ and $u_{2t-i}$ and $Q_{i,2\ell-i}$ connects $v_i$ and $v_{2\ell-i}$. 
We construct an $R$-parity-switcher 
  $$W=W(R,C_1,P_t, P_{1},\dots,\dots,P_{{t-1}},Q_{1},\dots,Q_{{\ell-1}}, Q_\ell, u_t,v_\ell).$$
Moreover, $|V(W)| \le 2+5(\ell+t-2) \le \frac{5}{12}\eta n +100$. 
This, combined with the fact that the edge $x_2y_2$ is disjoint with $V(W)$,  results in a contradiction with Claim~\ref{noR-parity}.

            When $|E(C_2) \cap E(R)|$ is odd, we label $P_{13}'$ as $(x_1=)u_1u_2\dots u_t \dots u_{2t}(=x_3)$ and $Q_{13}'$ as $(y_3=)v_1 v_2 \dots v_\ell \dots v_{2\ell}(=y_1)$ for some $t,\ell \in \{1,2,3\}$. It follows that $|E(C_1)|=2t+2\ell$.
       Similarly, we obtain $t+\ell -1$ disjoint paths of length at most $4$: $P_{1},\dots,P_{i},\dots,P_{{t-1}}, P_{x_3,x_4}$ and $Q_{1},\dots,Q_{i},\dots,Q_{{\ell-1}},Q_{y_4y_1}$ such that
            $P_{i}$ connects $u_i$ and $u_{2t-i}$, $P_{x_3,x_4}$ connects $x_3$ and $x_4$, $Q_{i}$ connects $v_i$ and $v_{2\ell-i}$, and $Q_{y_4y_1}$ connects $y_1$ and $y_4$. Let $P_{x_4y_1}=Q_{y_4y_1}+x_4y_4+P_{x_3x_4}$. Then $|E(P_{x_4y_1})|\le 9$.
Set $P_t=u_t$ and $Q_\ell=v_\ell$. We construct a $R$-parity-switcher 
            $$W=W(R,C_1,P_t, P_{1},\dots,P_{{t-1}},P_{x_4y_1},Q_{1},\dots,Q_{{\ell-1}}, Q_\ell, u_t,v_\ell).$$
 Additionally, $|V(W)| \le 2+5(\ell+t-1)+9 \le \frac{5}{12}\eta n +100$. 
 This, combined with the fact that the edge $x_2y_2$ is disjoint from $V(W)$, again leads to a contradiction to Claim~\ref{noR-parity}.

\end{proof}
The proof of this lemma is completed.

\end{proof}

\begin{lem} \label{bipartite}
    Let $\frac{1}{n} \ll \alpha  \ll \eta$ and let $G$ be an $n$-vertex graph with a vertex partition $V(G) = X \cup Y \cup Z$ satisfying: 
    \begin{equation*}
     \left(\frac{1}{2}- \alpha\right)n \le |Y| \le |X| \le \left(\frac{1}{2}+ \alpha\right)n, \,\,
        |Z| \le \alpha n,
    \end{equation*}
    and the following degree conditions:
    \begin{equation*}
        \Delta(G[X]) \le  \eta n,
        \min\{\delta_{X}(Y), \delta_{Y}(X)\} \ge  \left(\frac{1}{2}- 2\eta\right)n,\,\,
        \min\{\delta_{Y}(Z), \delta_{X}(Z)\} \ge \frac34\eta n.
    \end{equation*}
If
    \begin{equation} \label{-1}
         |X|-|Y|-|Z| \le f(X)-1
    \end{equation}
holds, where $$f(X) = \max\{e(\mathcal{P}): \mathcal{P} \text{ is  a UDP in } G[X]\},$$
   then
     $\mathcal{C}_{n}(G)=\mathcal{C}(G)$.

\end{lem}

\noindent\textbf{Remark.} The necessity of inequality (\ref{-1}) is demonstrated by the non-Hamilton-generated graph $G_1$ in Construction A. While this graph satisfies all  requirements specified in Lemma~\ref{bipartite},
with the exception that inequality (\ref{-1}) is replaced by the contradictory inequality $1=|X|-|Y|\le f(X)-1=0$, where $Z=\emptyset$ for $G_1$.

\begin{proof}
    Suppose that $\mathcal{C}(G) \neq \mathcal{C}_n(G)$. Again, by Lemma \ref{R-lem}, there exists a subgraph $R$ of $G$ satisfying conditions (a), (b), and (c). From (c), we have 
$$
d_R(x) \ge \frac{1}{2} d_G(x) \quad \text{for any } x \in V(G).
$$
Let $|X| - |Y| - |Z| = t$. Then we have 
$$|t| = |2|X| - n| \le 2\alpha n, \text{ and $t$ is odd.}$$
For any $x \in X$, we have
\begin{align}
    |N_R(x) \cap Y| &= d_R(x) - |N_R(x) \cap (X \cup Z)| \nonumber \\
    &\ge \frac{1}{2} d_G(x) - \eta n - \alpha n \nonumber \\
    &\ge \left(\frac{1}{4} - 3\eta \right)n \label{deltaYXn} \\
    &\ge \frac{\frac{1}{4} - 3\eta}{\frac{1}{2} + \alpha}|Y| \nonumber \\
    &\ge \left(\frac{1}{2} - 7\eta \right)|Y|, \label{deltaYX}
\end{align}
where the first inequality holds because $|N_R(x) \cap X| \le \Delta(G[X]) \le \eta n$ and $|N_R(x) \cap Z| \le |Z| \le \alpha n$.

    \begin{claim}
        $G[X,Y\cup Z]$ is $(\frac12\eta n,4)$-connected.
    \end{claim}
    \begin{proof}
        For an arbitrary subset $U \subseteq V(G)$ of size at most $\frac{1}{2} \eta n$, let $X' = X \setminus U$, $Y' = Y \setminus U$, $Z' = Z \setminus U$, and $G' = G[X', Y' \cup Z']$. Consider two vertices $v_1, v_2$ in $G'$. Since $\delta_X(Y), \delta_X(Z) \geq \frac{3}{4}\eta n > |U| + 2$, we can choose two additional vertices $v_1', v_2'$ in $V(G')$ such that for $i \in \{1, 2\}$:
$$
v_i' \in N_{G'}(v_i) \subseteq X \quad \text{if } v_i \in Y \cup Z, \quad \text{and } v_i' = v_i \quad \text{if } v_i \in X.
$$
Since $\delta_Y(X) \geq \left(\frac{1}{2} - 2\eta\right)n$, both $v_1'$ and $v_2'$ have at least 
$$
\left(\frac{1}{2} - 2\eta\right)n - |U| - 2 \geq \frac{\frac{1}{2} - 3\eta}{\frac{1}{2} - \alpha}|Y'| \geq (1-6\eta)|Y'|
$$
neighbors in $Y' \setminus \{v_1, v_2\}$. Thus, we can find at least one common neighbor of $v_1'$ and $v_2'$ in $Y' \setminus \{v_1, v_2\}$. This implies that 
$$
\mathrm{dist}(v_1, v_2) \leq 4,
$$ 
and the proof is completed.
        %$\delta(G') \ge \min\{\delta_X(Y),\delta_X(Z),\delta_Y(X)\} \ge \frac{3}{4}\eta n$
\end{proof}
    
    \begin{claim} \label{balanced-cycle}
        There does not exist a cycle $C$ in $G$ satisfying the following:
        \begin{enumerate}
        \item[(*)] $
        |X\cap V(C)|-|Y\cap V(C)|-|Z\cap V(C)|=0,
        $
        \item[(**)] $
        f(X\setminus V(C)) \ge |X\setminus V(C)|-|Y\setminus V(C)|-|Z\setminus V(C)|=t,
        $
        \item[(***)] $|E(C)|\le \frac17\eta n$ is even, and $|E(C)\cap E(R)|$ is odd.
        \end{enumerate}
    \end{claim}

    \begin{proof}
      If not, suppose $C=(v_1,v_2,...,v_{h},v_{h+1},...,v_{2h},v_1)$ is such a cycle, where $2h \le \frac17\eta n$.
        By (*),  we deduce that
        $$
        |X\setminus V(C)|-|Y\setminus V(C)|-|Z\setminus V(C)|=|X|-|Y|-|Z|=t.
        $$
     Choose a UDP $\mathcal{P}_0$ in $G[X\setminus V(C)]$ with $\max\{t,0\}$ edges, and two vertices $v_1' \in X\setminus \big(V(C) \cup V(\mathcal{P}_0)\big)$ and $v_{h+1}' \in Y\setminus \big(V(C) \cup V(\mathcal{P}_0)\big)$.
        %From $\delta_X(Y),\delta_Y(X) \ge (\frac12-2\eta)n$, which implies $e_G(X,Y) \ge (\frac{1}{4}-\eta)n^2$, we deduce that the bipartite graph $G[X, Y]$ is $(\sqrt{\eta}n, 4)$-connected by Lemma~\ref{bidiameter}~(2).
       Note that
       \begin{align}\label{V(P0)}
           |V(\mathcal{P}_0)|\le 2|t| \le 4\alpha n\le \frac{1}{14}{\eta}n.
       \end{align}
    We obtain that $G[X, Y\cup Z]- V(\mathcal{P}_0)$ is $(\frac{3}{7}{\eta}n, 4)$-connected. 
       
Note that $h+1 \le \frac{1}{14} \eta n + 1 \le \frac{3}{35} \eta n = \frac{m}{d+1}$.   Apply Lemma~\ref{bidiameter} to $G[X, Y \cup Z] - V(\mathcal{P}_0)$ for  $h + 1$ pairs of vertices 
$$
\{(v_1, v_1'), (v_{h+1}, v_{h+1}')\} \cup \left( \bigcup_{2 \le i \le h} \{(v_i, v_{2h+2-i})\} \right),
$$
We obtain, in $G[X, Y \cup Z] - V(\mathcal{P}_0)$, $h+1$ disjoint paths $P_1, P_2, \dots, P_h, P_{h+1}$ of length at most $4$ such that $P_i$ connects $v_i$ to $v_{2h+2-i}$ for $2 \le i \le h$,  and $P_i$ connects $v_i$ to $v'_{i}$ for $i \in \{1, h+1\}$.
Therefore, we construct an $R$-parity-switcher 
$$W = W(R, C, P_1, P_2,  \dots, P_h,P_{h+1}, v'_1, v'_{h+1}).$$
Let $$V' = \{v_1', v_{h+1}'\} \cup \left( V(G) \setminus \bigcup_{i=1}^{h+1} V(P_i) \right).$$
Denote $X' = V' \cap X$, $Y' = V' \cap Y$, and $Z' = V' \cap Z$.

We propose that $|X'| - |Y'| - |Z'| = t$. Consider the subgraph $H$ of $G[X, Y \cup Z]$ such that $V(H) = \bigcup_{i=1}^{h+1} V(P_i)$ and $E(H) = \bigcup_{i=1}^{h+1} E(P_i)$. Note that $H$ is a bipartite graph and $V(C) \subseteq V(H)$. It is sufficient to show that
\begin{align}\label{H:X,YZ}
    |V(H) \cap X| = |V(H) \cap (Y \cup Z)|
\end{align}
since $v_1' \in X$ and $v_{h+1}' \in Y$. Observe that $d_H(u) = 1$ when $u \in V(C) \cup \{v_1', v_{h+1}'\}$, and $d_H(u) = 2$ otherwise. Since $|(V(C) \cup \{v_1', v_{h+1}'\}) \cap X| = |(V(C) \cup \{v_1', v_{h+1}'\}) \cap (Y \cup Z)|$, we deduce (\ref{H:X,YZ}) by the equality
$$
\sum_{u \in V(H) \cap X} d_H(u) = \sum_{u \in V(H) \cap (Y \cup Z)} d_H(u).
$$

Denote the endpoints of paths in $\mathcal{P}_0$ as $w_i^1, w_i^2$ for $i = 1, 2, \dots, |\mathcal{P}_0|$, and let 
$$ r = |Z'| + \min\left\{ 0, \frac{t - 1}{2} \right\}. $$
Label $Z' = \{ z_1, z_2, \dots, z_{|Z'|} \}$. By 
$$ \delta_X(Z), \delta_Y(Z) \ge \frac{3}{4} \eta n \ge 6 \alpha n + 1 \ge 2 |Z| + |V(\mathcal{P}_0)| + 1, $$ 
we can    iteratively and greedily find pairwise disjoint neighbors $u_i^1, u_i^2$ in $X' \setminus \left( \{v_1'\} \cup V(\mathcal{P}_0) \right)$ for each $z_i \in Z'$ with $i \in [r]$.
If $r < |Z'|$, we further assign $u_i^1, u_i^2 \in Y'$ for $i \in \{r + 1, r + 2, \dots, |Z'| - 1\}$, and 
$ u_{|Z'|}^1 \in X' \setminus \left( \{v_1'\} \cup V(\mathcal{P}_0) \right), u_{|Z'|}^2 \in Y'$.

Let $X'' = X' \setminus \text{In}(\mathcal{P}_0)$. Note that
\begin{eqnarray*}
    |X''| &\ge& |X\setminus (\cup_{1\le i \le h+1}V(P_i) \cup V(\mathcal{P}_0))|\\ 
         &\ge& \left(\frac12-\alpha\right)n-\left(5\cdot\left(\frac{1}{14}\eta n+1\right)+\frac{1}{14}\eta n\right) \\ 
         &\ge& \frac12\left(1-\eta\right)n
\end{eqnarray*}
and $ |Y'| \ge |Y\setminus (\cup_{1\le i \le h+1}V(P_i))| \ge \frac12(1-\eta)n$. We have $|V(G) \setminus (X''\cup Y')| \le \eta n$.

Now we want to apply Lemma \ref{Fixpair} (2) to the bipartite graph $G'=G[X'',Y']$.
%for $\varepsilon_1 = 3\eta, \varepsilon_2=6\alpha$, 
Note that $\delta(G')\ge \min\{\delta_{X''}(Y'),\delta_{Y'}(X'')\} \ge (\frac12-2\eta)n-|V(G)\setminus (X'' \cup Y')| \ge (\frac12-3\eta)n$. 
We choose $\ell = 1+|Z'|+|\mathcal{P}_0| \le 1+\alpha n + 4\alpha n \le 6\alpha  n$ pairs (by (\ref{V(P0)})) of vertices: 
$$
\{(v_1',v_{h+1}')\} \cup \bigcup_{i=1}^{|Z'|}\{(u_i^1,u_i^2)\} \cup \bigcup_{i=1}^{|\mathcal{P}_0|}\{(w_i^1,w_i^2)\}.
$$
Label the set of vertices in above pairs as $S$. Then 
$$ S\cap X''=\{v_1'\}\cup\bigcup_{i=1}^{r}\{u_i^1,u_i^2\}\cup\bigcup_{i=1}^{|\mathcal{P}_0|}\{w_i^1,w_i^2\}$$
or $$ S\cap X''=\{v_1'\}\cup\bigcup_{i=1}^{r}\{u_i^1,u_i^2\}\cup\bigcup_{i=1}^{|\mathcal{P}_0|}\{w_i^1,w_i^2\}\cup\{u_{|Z'|}^1\} \text{ (when $r < |Z'|$)},$$
%\begin{itemize}
 %       \item $v_1'$;
  %      \item $u_1^1,u_1^2,u_2^1,u_2^2,...,u_r^1,u_r^2$;
   %     \item $w_1^1,w_2^1,w_2^1,w_2^2,...,w_{|\mathcal{P}_0|}^1,w_{|\mathcal{P}_0|}^2$;
    %    \item $u_{|Z'|}^1$ if $r < |Z'| $,
    %\end{itemize}
and $$S \cap Y'=\{v_{h+1}'\}\cup\bigcup_{i=r+1}^{|Z'|-1}\{u_i^1,u_i^2\} $$ 
or $$S \cap Y'=\{v_{h+1}'\}\cup\bigcup_{i=r+1}^{|Z'|-1}\{u_i^1,u_i^2\}\cup\{u_{|Z'|}^2\} \text{ (when $r < |Z'|$)}. $$
%\begin{itemize}
 %       \item $v_{h+1}'$;
  %      \item $u_{r+1}^1,u_{r+1}^1,...,u_{|Z'|-1}^1,u_{|Z'|-1}^2,u_{|Z'|}^2$ if $r < |Z'|$.
   % \end{itemize}
Thus by $r=|Z'|+\min\{0,\frac{t-1}{2}\}$, we know
\begin{eqnarray*}
    &|S\cap X''|&=1 + \min\{2|Z'|,2r+1\}+2|\mathcal{P}_0|=1+2|Z'|+\min\{0,t\}+2|\mathcal{P}_0|,\\
    &|S \cap Y'|&=1+\max\{0,2|Z'|-2r-1\}=1+\max\{0,-t\}.
\end{eqnarray*}
Therefore, 
    \begin{eqnarray*}
    |X''|-|Y'|-\frac12(|S\cap X''|-|S \cap Y'|) &=& |Z'|+t-|\text{In}(\mathcal{P}_0)|-(|Z'|+\min\{0,t\}+|\mathcal{P}_0|) \\
    &=& \max\{t,0\}-|E(\mathcal{P}_0)|\\ &=& 0.
    \end{eqnarray*}
The condition (\ref{balanced}) of Lemma~\ref{Fixpair} (2) is satisfied.
Let \[
E_0 = \{ v_1' v_{h+1}' \} \cup \bigcup_{i=1}^{|Z'|} \{ u_i^1 u_i^2 \} \cup \bigcup_{i=1}^{|\mathcal{P}_0|} \{ w_i^1 w_i^2 \}.
\]
Applying  Lemma~\ref{Fixpair} (2) to the bipartite graph $G'=G[X'',Y']$ with the vertex pairs of $E_0$, and parameters $\varepsilon_1= 3\eta$, $\varepsilon_2= 6\alpha$, we obtain a Hamiltonian cycle $C$ containing \( E_0 \) in the graph \( G[X'', Y'] + E_0 \).
Therefore we have a $v_1'$-$v_{h+1}'$-path $P_{v_1'v_{h+1}'}$   that covers exactly $V'$ from $C$ by replacing each \( u_i^1 u_i^2 \) by the path \( u_i^1 z_i u_i^2 \) and each \( w_i^1 w_i^2 \) by its corresponding path in \( \mathcal{P}_0 \), and deleting the edge \( v_1' v_{h+1}' \).
This results in a contradiction to Corollary \ref{parity}.
    \end{proof}

    \begin{claim} \label{same-parity}
        Let $\mathcal{P}_0$ be a UDP in $G[X]$ with $\max\{t,0\}$ edges. Then there do not exist vertices $x_1,x_2 \in X\setminus V(\mathcal{P}_0)$ and $y_1\in Y \cup Z$ such that $x_1y_1 \in E(R)$ and $x_2y_1 \in E(G) \setminus E(R)$.
    \end{claim}
    \begin{proof}
        Suppose this is the case. Let $X'=X\setminus V(\mathcal{P}_0)$. 
Since for any $S \subseteq X'$,  $\mathcal{P}_0$ is still an UPD in $G[X\setminus S]$, we have $f(X\setminus S) \ge |E(\mathcal{P}_0)| = \max\{t,0\}$.
        %that $x_1,x_2 \in X\setminus V(P_0),y\in Y \cup Z$ such that $x_1y \in E(R), x_2y \in E(G) \setminus E(R)$. 
        If $N_R(x_1) \cap N_R(x_2) \cap Y \ne \emptyset$, choose $y_0$ in it.
        Then we construct an even cycle $C=(x_1,y_1,x_2,y_0,x_1)$ satisfying (*), (**), (***) in Claim~\ref{balanced-cycle}, a contradiction. Hence $N_R(x_1) \cap N_R(x_2) \cap Y  = \emptyset$. Therefore, by (\ref{deltaYX}), we have $|(N_R(x_1) \cup N_R(x_2)) \cap Y | \ge (1-14\eta)|Y|$.
        %, and recall that ${\color{blue}|N_R(x_1) \cap Y|,|N_R(x_2) \cap Y| \ge (1/2-2\eta)|Y|}$
If there exists $x_3 \in X'\setminus\{x_1,x_2\}$ such that 
        $$\min\{|N_R(x_1) \cap N_R(x_3) \cap Y|, |N_R(x_2) \cap N_R(x_3) \cap Y|\} \ge 3,$$
        then there are two vertices $y_2,y_3 $ such that
        $y_i \in (N_R(x_3)\cap N_R(x_{i-1}) \cap Y)\setminus \{y_1\}$ for $i \in \{2,3\}$.
        Therefore, $(x_1,y_1,x_2,y_3,x_3,y_2,x_1)$ forms a forbidden cycle in Claim \ref{balanced-cycle}. 
This implies that for any $x \in X'\setminus \{x_1,x_2\}$, we have either
        $|N_R(x) \cap N_R(x_1) \cap Y| \le 2$
        or $|N_R(x) \cap N_R(x_2) \cap Y| \le 2$.
        Moreover, recall that $|N_R(x)\cap Y| \ge (1/2-7 \eta)|Y|$ for any $x \in X$ by (\ref{deltaYX}), we know that exactly one of these two cases holds. 
        
        For \( i \in \{1,2\} \), let \( Y_i = N_R(x_i) \cap Y \). Then 
        \begin{align*}
             |Y_i| =|N_R(x_i) \cap Y| \ge \left(\frac{1}{4} - 4\eta \right)n.
          \end{align*}
        Partition $X'$ into $X_1$ and $X_2$, where
    \begin{align*}
        X_1=\{x \in X': |N_R(x) \cap Y_2| \le 2 \} \quad \text{and} \quad X_2=\{x \in X': |N_R(x) \cap Y_1| \le 2 \}.
    \end{align*}
   Let \( V_i = X_i \cup (N_R(x_i) \cap Y) \) for \( i \in \{1,2\} \). 
          Note that \( Y_1 \) and \( Y_2 \) are disjoint by their definitions. Hence $V_1$ and $V_2$ are disjoint too. 
          
          We now give an upper bound of \( e_R(V_1, V_2) \) as:
\begin{align*}
    e_R(V_1, V_2) &\le e_G(X) + |Y_1| \cdot |Y_2| + e_R(X_1, Y_2)+e_R(X_2,Y_1) \\
    &\le \frac{1}{2}\eta n^2 + \frac{1}{4}|Y|^2 + 2|X| \\
    &\le \left(\frac{1}{16} + \eta\right)n^2,
\end{align*}
where \( e_R(X_1, Y_2)+e_R(X_2,Y_1) \le 2|X| \) follows from the properties of \( X_1 \) and \(X_2\).

Next, we provide a lower bound for \( e_G(V_1, V_2) \):
\begin{align*}
    e_G(V_1, V_2) &\ge e_{G \setminus R}(X_1, Y_2) + e_{G \setminus R}(X_2, Y_1) + e_R(V_1, V_2) \\
    &\ge \left(|X_1| + |X_2|\right)\left(\frac{1}{4} - 7\eta\right)n + e_R(V_1, V_2) \\
    &\ge \left(\frac{1}{8} - 4\eta\right)n^2 + e_R(V_1, V_2).
\end{align*}
Here, we use the fact that for any $i \in \{1,2\}$ and $x \in X_i$, 
\begin{align*}
    |N_{G\setminus R}(x) \cap Y_{3-i}|
    &\ge |N_G(x) \cap Y_{3-i}| - 2\\ 
    &\ge \left(\frac12-2\eta\right)n-2-|Y \setminus Y_{3-i}|\\
    &\ge \left(\frac12-2\eta\right)n-2-\left(\frac12+\alpha\right) n+\left(\frac{1}{4} - 4\eta \right)n\\
    &\ge \left(\frac{1}{4} - 7\eta \right)n.
\end{align*}
Combining the above two bounds, we obtain:
\[
\left(\frac{1}{8} - 4\eta\right)n^2 \le e_G(V_1, V_2) - e_R(V_1, V_2) \le e_R(V_1, V_2) \le \left(\frac{1}{16} + \eta\right)n^2,
\]
which leads to a contradiction.
%to the property (c) of $R$.
    \end{proof}

    \begin{claim} \label{E(X,Y)}
          Suppose $\mathcal{P}_0$ is a UDP in $G[X]$ such that $|E(\mathcal{P}_0)|=\max\{t+1,0\}$. Then there exists a partition $Y \cup Z=V_1 \cup V_2$ such that the followings hold:
            \begin{itemize}
                \item[(i)] $|V_1|\le 4\eta n$;
                \item[(ii)] $E_G(X\setminus \text{In}(\mathcal{P}_0),V_1)\cap E(R)=\emptyset$ and $E_G(X\setminus \text{In}(\mathcal{P}_0),V_2) \subseteq E(R)$;
                \item[(iii)] $\big(E_G(V_1)\cup E_G(V_2)\big) \cap E(R) = \emptyset$ and $E_G(V_1,V_2) \subseteq E(R).$
            \end{itemize}

            %$$E_G(X\setminus In(\mathcal{P}_0),V_1)\cap E(R)=\emptyset,E_G(X\setminus In(\mathcal{P}_0),V_2) \subseteq E(R)$$ and $$\big(E_G(V_1)\cup E_G(V_2)\big) \cap E(R) = \emptyset, E_G(V_1,V_2) \subseteq E(R), |V_1|\le 4\eta n.$$
        \end{claim}

        \begin{proof}
            For any 
            $x \in \text{End}(\mathcal{P}_0)$, denote $\mathcal{P}_x=\mathcal{P}_0 \setminus \{x\}$. Then $\mathcal{P}_x$ is a UDP with $|E(\mathcal{P}_x)|=\max\{0,t\}$. By Claim \ref{same-parity}, 
            we have for any $v \in Y \cup Z$, either
            $$E_G(v,X\setminus V(\mathcal{P}_x)) \subseteq E(R) \quad \text{or} \quad E_R(v,X\setminus V(\mathcal{P}_x)) = \emptyset.$$
            Since $x$ is chosen randomly in $\text{End}(\mathcal{P}_0)$,
            we can partition $Y \cup Z$ into $V_1\cup V_2$, where
            \begin{align*}
               V_1 &= \{v \in Y \cup Z : E_R(v, X \setminus \text{In}(\mathcal{P}_0)) = \emptyset\}, \\
               V_2 &= \{v \in Y \cup Z : E_G(v, X \setminus \text{In}(\mathcal{P}_0)) \subseteq E(R)\}.
            \end{align*}
            This means that $E_R(V_1,X\setminus \text{In}(\mathcal{P}_0))=\emptyset $ and $E_G(V_2,X\setminus \text{In}(\mathcal{P}_0)) \subseteq E(R)$, i.e. the statement (ii) holds.
            By (\ref{deltaYX}), we have $|V_2| \ge (1/2-7\eta)|Y|$.

            Next, we prove that 
            $\big(E_G(V_1)\cup E_G(V_2)\big) \cap E(R) = \emptyset$ and $E_G(V_1,V_2) \subseteq E(R).$
            By (\ref{V(P0)}),  we have $|V(\mathcal{P}_0)|\le 2|t|+2 \le 4\alpha n +2 <\min\{\delta_X(Y), \delta_X(Z)\}-2$. For arbitrary $f_0 = ab \in E(G[Y \cup Z])$, we can choose two vertices $a' \in N_G(a)\cap (X\setminus V(\mathcal{P}_0))$ and $b' \in N_G(b)\cap (X\setminus V(\mathcal{P}_0))$. 
            Let End$(\mathcal{P}_0)=\bigcup_{i=1}^{|\mathcal{P}_0|}\{w_i^1, w_i^2\}$.
           % the endpoints of paths in $\mathcal{P}_0$ as $w_i^1, w_i^2$ for $i = 1, 2, \dots, |\mathcal{P}_0|$.
            %Then we can find disjoint neighbors, denoted as $\bigcup_{i=1}^{|\mathcal{P}_0|} \{w_i^1,w_i^2\}$, in $V_2 \setminus (Z\cup \{a,b\})$ of $\bigcup_{i=1}^{|P_0|} \{v_i^1,v_i^2\}$, so $\{v_i^jw_i^j\} \subseteq E(R)$ for all $i,j$ by the definition of $V_1,V_2$.
            Denote $Z'=Z\setminus \{a,b\}=\{z_1,z_2,...,z_{|Z'|}\}$ and \[r=|Z'|+\min \left\{0, \frac{t+1}{2}\right\}.\]
           Since $\delta_X(Z),\delta_Y(Z) \ge \frac34 \eta n$ and $|Z| \le \alpha n$, we can iteratively and greedily find pairwise disjoint neighbors $u_i^1,u_i^2$ for each  $z_i\in Z'$ such that for all $1\le i \le |Z'|$, the followings hold:
            \begin{itemize}
                \item $\{z_iu_i^1,z_iu_i^2 \}\subseteq E(R)$ or $\{z_iu_i^1,z_iu_i^2 \}\subseteq E(G)\setminus E(R)$;
                \item $\{u_i^1,u_i^2\} \subseteq X \setminus (V(\mathcal{P}_0)\cup \{a',b'\})$ for $i \in \{1,2,...,r\}$;
                \item $\{u_i^1,u_i^2\} \subseteq V_1 \cap Y$ or $V_2 \cap Y$ for $i \in \{r+1,r+2,...,|Z'|\}$.
            \end{itemize}
    Now let $X'=X\setminus \text{In}(\mathcal{P}_0)$ and $Y'=Y\setminus \{a,b\}$. Then 
    \begin{align*}
        \delta_{Y'}(X') \ge \delta_Y(X)-2 \ge \left(\frac12-3\eta\right)n \text{ and } \delta_{X'}(Y') \ge \delta_X(Y)- \max\{0,t\} \ge \left(\frac12-3\eta\right)n.
    \end{align*}
     Choose $1+|Z'|+|\mathcal{P}_0|(\le 3\alpha n)$ pairs of vertices
            $$
            \{(a',b')\} \cup \bigcup_{i=1}^{|Z'|}\{(u_i^1,u_i^2)\} \cup \bigcup_{i=1}^{|\mathcal{P}_0|} \{(w_i^1,w_i^2)\}.
            $$
            Let $S= \{a',b'\} \cup \left(\bigcup_{i=1}^{|Z'|}\{u_i^1,u_i^2\}\right) \cup \left(\bigcup_{i=1}^{|\mathcal{P}_0|} \{w_i^1,w_i^2\}\right)$  and 
            $$
            E_0 = \{ aa',ab,bb' \} \cup \bigcup_{i=1}^{|Z'|} \{ u_i^1 u_i^2 \} \cup \bigcup_{i=1}^{|\mathcal{P}_0|} \{ w_i^1 w_i^2 \}.
            $$
    Since  
            \begin{eqnarray*}
              & & |X'|-|Y'|-\frac12(|S\cap X'|-|S \cap Y'|)\\ 
              &\quad & =|Z'|+t+|\mathcal{P}_0|+2-\max\{t+1,0\}-(1+|\mathcal{P}_0|+|Z'|+\min\{t+1,0\})\\
               &\quad &= t+2-\max\{t+1,0\}-1-\min\{t+1,0\}\\
               &\quad &= 0,
            \end{eqnarray*}
      the condition (\ref{balanced}) is satisfied.
     Apply Lemma \ref{Fixpair}~(2) to $G[X',Y']$ with parameters $\varepsilon_1= 3\eta$, $\varepsilon_2= 3\alpha$,
     we have a Hamilton cycle $C'$ containing \( E_0 \) in the graph \( G[X', Y'] + E_0 \).
   Therefore, we have    a Hamilton cycle $C$  obtained from $C'$ by replacing each \( u_i^1 u_i^2 \) by the path \( u_i^1 z_i u_i^2 \) and each \( w_i^1 w_i^2 \) by its corresponding path in \( \mathcal{P}_0 \), and $a'b'$ by the path $a'abb'$.
     
            From the construction of $C$, we have $E_G(C) \cap E_G(X', Y')=E_G(C) \setminus E_0$. 
            We assert that %its intersection with $E(R)$ consists an even number of edges.
            $|(E_G(C) \setminus E_0) \cap E(R)|$ is even. 
           % This holds because:
%Let $H$ be the subgraph induced by the edge set $E_G(C) \setminus E_0$. 
For any vertex  $y\in Y' \setminus \bigcup_{i=r+1}^{|Z'|} \{u_i^1, u_i^2\}$,  
there are precisely two edges incident to $y$ in $E_G(C) \setminus E_0$. These two edges are either both contained in $E(R)$ or both in $E(G) \setminus E(R)$, as a result of the properties of $V_1$ and $V_2$.
For $u_i^1$ and $u_i^2$ with $i \in \{r+1, \dots, |Z'|\}$, there are precisely one edge incident to $u_i^j$ in $E_G(C) \setminus E_0$ for $j=1,2$. Moreover, $u_i^1$ and $u_i^2$ are either both in $V_1$ or both in $V_2$.  
Therefore, these two edges incident with $u_i^1$ and $u_i^2$ are either both in $E(R)$ or both in $E(G) \setminus E(R)$. 
Hence we have $|(E_G(C) \setminus E_0) \cap E(R)|$ is even.
%For any $i \in \{r+1, \dots, |Z'|\}$, we know that both $u_i^1$ and $u_i^2$ have a degree of $1$, {\color{red}and are in either $V_1$ or $V_2$}, 
%so these two edges incident with them are both in either $E(R)$ or $E(G) \setminus E(R)$. 
%Hence, $|(E_G(C) \setminus E_0) \cap E(R)|$ is even.

%Additionally, for any \( i \in \{r+1, \dots, |Z'|\} \), the edges \( \{z_i u_i^1, z_i u_i^2\} \) are contained either in \( E(R) \) or in \( E(G) \setminus E(R) \). 

Combining the above assertion with Lemma~\ref{R-lem}~(b), we conclude that the number of the remaining edges in $E_G(C) \cap E(R)$ is odd. Recall that $\{z_i u_i^1, z_i u_i^2\}$ always contributes an even number of edges in $E(R)$. We deduce that
$$
\big|\left(\{aa', ab, bb'\} \cup E_G(\mathcal{P}_0)\right) \cap E(R)\big|
$$
must be even. Since $\mathcal{P}_0$ is fixed, the parity of $|\{aa', ab, bb'\} \cap E(R)|$ is independent of the choice of these four vertices. 
By (ii), for any $x \in \{a,b\}$, we have 
\begin{align}\label{a'b'}
    xx' \in E(G)\setminus E(R) \text{ if } x\in V_1, \text{ and } xx' \in E(R) \text{ if } x\in V_2.
\end{align}
If $|\{aa', ab, bb'\} \cap E(R)|$ is even, by (\ref{a'b'}) we know that
$ab \in  E(G)\setminus E(R)$ if $a,b$ are both in $V_1$ or both in $V_2$, and 
$ab \in E(R)$ if one of $a,b$ is in $V_1$ and another is in $V_2$. Hence,
$$\big(E_G(V_1) \cup E_G(V_2)\big) \cap E(R) = \emptyset \quad \text{and} \quad E_G(V_1, V_2) \subseteq E(R).$$
And similarly, if $|\{aa', ab, bb'\} \cap E(R)|$ is odd, 
$$E_G(V_1), E_G(V_2) \subseteq E(R) \quad \text{and} \quad E_G(V_1, V_2) \cap E(R) = \emptyset.$$

%{\color{red}Thus, we can deduce the following:
%\begin{align*}
 %   \text{either} \quad &\big(E_G(V_1) \cup E_G(V_2)\big) \cap E(R) = \emptyset \quad \text{and} \quad E_G(V_1, V_2) \subseteq E(R), \\
  %  \text{or} \quad &E_G(V_1), E_G(V_2) \subseteq E(R) \quad \text{and} \quad E_G(V_1, V_2) \cap E(R) = \emptyset.
%\end{align*}
%}
Now we prove that the second case can not hold. Otherwise, combined with (ii), we would have 
$$E_G\big(V_1, V_2 \cup X \setminus \text{In}(\mathcal{P}_0)\big) \cap E(R) = \emptyset,$$ 
which implies
\[
e_R(V_1, \overline{V_1}) \leq e_G(V_1, \text{In}(\mathcal{P}_0)) \leq |V_1| \cdot |\text{In}(\mathcal{P}_0)| \leq \max\{0,t\}\cdot |V_1| \leq 2\alpha n|V_1| \leq \frac{\eta n}{3}|V_1|.
\]
However, from the property (c) of $R$, we also have
\[
2e_R(V_1, \overline{V_1}) \geq e_G(V_1, \overline{V_1}) \geq e_G(V_1, X) \geq \delta_X(Z)|V_1| \geq \eta n|V_1|.
\]
These two bounds can hold simultaneously only if \( |V_1| = 0 \), which would imply \( E(G) = E(R) \), a contradiction.
The statement (iii) holds.

Finally, we prove that \( |V_1| \leq 4\eta n \). To do this, it suffices to show \( |V_1 \setminus Z| \leq 3\eta n \) since \( |Z| \leq \alpha n \). 
Combine with (ii) and (iii), we have 
\begin{align*}
    e_R(V_1 \setminus Z, \overline{V_1 \setminus Z})
    &= e_R(V_1 \setminus Z, \text{In}(\mathcal{P}_0) \cup V_2) \\
    &\leq |V_1 \setminus Z| \big(|\text{In}(\mathcal{P}_0)| + |V_2|\big) \\
    &\leq |V_1 \setminus Z| \big(\max\{t, 0\} + |Y| + |Z| - |V_1 \setminus Z| \big) \\
    &= |V_1 \setminus Z| \big(\max\{|X|, |Y| + |Z|\} - |V_1 \setminus Z| \big).
\end{align*}
Moreover, 
\begin{align*}
    e_G(V_1 \setminus Z, \overline{V_1 \setminus Z}) 
    &\geq e_R(V_1 \setminus Z, \text{In}(P_0) \cup V_2) + e_G(V_1 \setminus Z, X \setminus \text{In}(P_0)) \\
    &\geq e_R(V_1 \setminus Z, \overline{V_1 \setminus Z}) + |V_1 \setminus Z| \big(\delta_X(Y) - |\text{In}(P_0)| \big) \\
    &\geq e_R(V_1 \setminus Z, \overline{V_1 \setminus Z}) + (1 - 4\eta) |V_1 \setminus Z| \cdot |X|.
\end{align*}
Combining the above two inequalities and using \( e_G(V_1 \setminus Z, \overline{V_1 \setminus Z}) \leq 2e_R(V_1 \setminus Z, \overline{V_1 \setminus Z}) \), we obtain
\[
 (1 - 4\eta) |V_1 \setminus Z|\cdot |X| \leq |V_1 \setminus Z| \big(\max\{|X|, |Y| + |Z|\} - |V_1 \setminus Z| \big).
\]
Therefore, \( |V_1 \setminus Z| \leq 4\eta |X| + \max\{0,t\} \leq 3\eta n \), where the second inequality holds since $|X| \le (\frac12+\alpha)n$ and $|t| \le 2\alpha n$.
\end{proof}

Choose \( \mathcal{P}_0 \) as a UDP in \( G[X] \) such that \( |E(\mathcal{P}_0)| = \max\{t+1, 0\} \). Then, by Claim~\ref{E(X,Y)}, \begin{align}\label{X1X2}
\noindent{\text{we have a partition of $Y\cup Z=V_1\cup V_2$ satisfying (i), (ii), and (iii). }   }
\end{align}

\begin{claim} \label{X-path}
    For any UDP \( \mathcal{P} \) in \( G[X] \) such that \( |E(\mathcal{P})| = \max\{t, 1\} \), we have that \( |E(\mathcal{P}) \cap E(R)| \) is even and
    \begin{align*}
        E_G\big(X\setminus \text{In}(\mathcal{P}), V_1\big) \subseteq E(G) \setminus E(R)\quad \text{and} \quad E_G\big(X\setminus \text{In}(\mathcal{P}), V_2\big) \subseteq E(R).
    \end{align*}
\end{claim}

           \begin{proof}
Let \( I = V(\mathcal{P}_0) \setminus V(\mathcal{P}) \).
Label the endpoints of each path in \( \mathcal{P} \) as \( \{v_i^1, v_i^2\} \) for \( i \in \{1, 2, \dots, |\mathcal{P}|\} \), and label the vertices in \( I \) as 
\[
\{v_{|\mathcal{P}|+1}, v_{|\mathcal{P}|+2}, \dots, v_{|\mathcal{P}|+|I|}\}.
\]
For each \( i \in \{1, 2, \dots, |I|\} \), set \( v_{|\mathcal{P}|+i} = v_{|\mathcal{P}|+i}^1 = v_{|\mathcal{P}|+i}^2 \).
Recall that $|t|\le 2\alpha n$. Then we have
\begin{align*}%\label{P&P0}
    |\mathcal{P}|+|I| \le |V(\mathcal{P})|+|V(\mathcal{P}_0)| \le 2(|E(\mathcal{P})|+|E(\mathcal{P}_0)|) \le 4|t|+4 \le 9\alpha n.
\end{align*}

For any two vertices $x$ and $y$ satisfying $y \in N_G(v_1^1) \cap (Y\cup Z)$ and $x \in (N_G(y) \cap X)\setminus (I\cup V(\mathcal{P}))$. By combining Claim~\ref{E(X,Y)} (i), which follows $|V_1| \le 4\eta n$,
and (\ref{deltaYXn}), we have for any \( a \in X \) we have
\begin{align}\label{V2-Z}
    |N_R(a) \cap (V_2\setminus Z)|
    \ge \left( \frac{1}{4} - 3\eta \right)n-4\eta n-\alpha n
    \ge \left( \frac{1}{4} - 8\eta \right)n
    \ge 2(|\mathcal{P}|+|I|).
\end{align}
Therefore, we can greedily choose
the following disjoint pairs of vertices in $V_2$:
\[
\{(w_1^1, w_1^2), (w_2^1, w_2^2), \dots, (w_{|\mathcal{P}|+|I|}^1, w_{|\mathcal{P}|+|I|}^2)\},
\]
such that 
\begin{align}\label{subseteq}
    \{xw_1^1,v_1^2w_1^2\} \cup \{v_i^j w_i^j: 2 \le i \le |\mathcal{P}|+|I|, 1\le j \le 2\} \subseteq E(R).
\end{align}
Let $S_1=\{y\} \cup \bigcup_{i=1}^{|\mathcal{P}|+|I|}\{w_i^1, w_i^2\}$, 
$X'=X\setminus (I\cup V(\mathcal{P})\cup \{x\}),Y'=Y\setminus \{y\},Z'=Z\setminus\{y\}:=\{z_1,z_2,\dots,z_{|Z'|}\}$
and
\[
r = |Z'| + \min\left\{0, \frac{t-1}{2}\right\}.
\]
Since \( |Z| \leq \alpha n \) and \( \delta_X(Z), \delta_Y(Z) \geq \frac{3}{4} \eta n \geq 6|Z| \), we can greedily find pairwise disjoint pairs of vertices
\[
\{(u_1^1, u_1^2), (u_2^1, u_2^2), \dots, (u_{|Z'|}^1, u_{|Z'|}^2)\},
\]
such that the following conditions hold:
\begin{itemize}
    \item For  $i \in \{1, 2, \dots, r\},  \{u_i^1, u_i^2\} \subseteq X'$;
    \item For \( i \in \{r+1, r+2, \dots, |Z'|\} \), either \( \{u_i^1, u_i^2\} \subseteq V_1 \setminus \left(Z \cup S_1\right) \) or \( \{u_i^1, u_i^2\} \subseteq V_2 \setminus \left(Z \cup S_1\right) \);
    \item For all \( i \in \{1, 2, \dots, |Z'|\} \), either \( \{z_i u_i^1, z_i u_i^2\} \subseteq E(R) \) or \( \{z_i u_i^1, z_i u_i^2\} \subseteq E(G)\setminus E(R) \).
\end{itemize}
Note that $X'=X\setminus (V(\mathcal{P})\cup I \cup \{x\}) \subseteq X\setminus V(\mathcal{P}_0)\subseteq X\setminus \text{In}(\mathcal{P}_0)$.
Consider the bipartite graph $G':=G[X', Y']$.  Then 
$$\delta(G') \ge \delta_X(Y)-|V(\mathcal{P})\cup I|\ge \left(\frac12-3\eta\right)n\ge \left(\frac12-3\eta\right)|V(G')|.$$
Choose $\ell = |\mathcal{P}|+|I|+|Z| \ (\le 9\alpha n+\alpha n \le  10\alpha n)$ pairs of vertices
$$
\bigcup_{i=1}^{|\mathcal{P}|+|I|} \{(w_i^1,w_i^2)\} \cup \bigcup_{i=1}^{|Z'|} \{(u_i^1,u_i^2)\}.
$$
Label the set of vertices in above pairs as $S$, and let
$$
E_0=\bigcup_{i=1}^{|Z'|} \{ u_i^1 u_i^2 \} \cup \bigcup_{i=1}^{|\mathcal{P}|+|I|} \{ w_i^1 w_i^2 \}.
$$
Then 
\begin{align*}
& |X'| - |Y'| - \frac12(|S\cap X'|-|S\cap Y'|)\\ 
& = |X|-|V(\mathcal{P})|-|I|-1-|Y'|-(2r-|Z'|-|\mathcal{P}|-|I|)\\
& =|X|-|Y'\cup Z'|-1+2|Z'|-|E(\mathcal{P})|-2r\\
& =t+ 2|Z'|-\max\{t,1\}-2\left(|Z'|+\min \left\{0,\frac{t-1}{2}\right\}\right)\\
& =0.
\end{align*}
Hence (\ref{balanced}) of Lemma \ref{Fixpair} (2) holds.               
  Applying Lemma \ref{Fixpair} (2) to $G'$ with parameters $\varepsilon_1= 3\eta$, $\varepsilon_2= 10\alpha$, and $\ell$ pairs of vertices  
  $$ \bigcup_{i=1}^{|\mathcal{P}|+|I|} \{(w_i^1,w_i^2)\} \cup \bigcup_{i=1}^{|Z'|} \{(u_i^1,u_i^2)\}, $$          
we obtain a Hamiltonian cycle $C'$ containing \( E_0 \) in the graph \( G[X', Y'] + E_0 \). Therefore, we have a Hamilton cycle $C$ of $G$ obtained from $C'$  by replacing:
\begin{itemize}
    \item each \( u_i^1 u_i^2 \) by the path \( u_i^1 z_i u_i^2 \);
    \item $w_1^1 w_1^2$ by path $w_1^1xyv_1^1P_1v_1^2w_1^2$;
    \item each \( w_i^1 w_i^2 \) with $2 \le i \le |\mathcal{P}|$ by its corresponding path $w_i^1v_i^1P_iv_i^2w_i^2$ in \( \mathcal{P} \) ;
    \item and each \( w_i^1 w_i^2 \) with $ |\mathcal{P}|+1\le i \le |\mathcal{P}|+|I|$ by the path \( w_i^1 v_i w_i^2 \),
\end{itemize}
where $P_i$ is the path in $\mathcal{P}$ with endpoints $v_i^1$ and $v_i^2$.
We call the above process of construction of $C$ as an {\em HC Construction}.

We assert that the number $|(E(\mathcal{P})\cup \{v_1^1y,xy\})\cap E(R)|$ must be even. Note that
$$E_G(C)=E(\mathcal{P})\cup \{xy,yv_1^1\}\cup (E(C')\setminus E_0) \cup E_1,
$$
where $E_1=\{xw_1^1,v_1^2w_1^2\} \cup \bigcup_{i=1}^{|Z'|}\{u_i^1z_i, z_iu_i^2\}\cup\bigcup_{i=2}^{|\mathcal{P}|+|I|}\{v_i^1w_i^1, v_i^2w_i^2\}$ contains exactly an even number of edges in $E(R)$,
because we can combine with (\ref{subseteq}) and the fact that $\{z_iu_i^1, z_iu_i^2\}$ is either in $E(R)$ or in $E(G) \setminus E(R)$ for any $i \in [|Z'|]$.
Hence, to prove this it is sufficient to prove that $|(E_G(C') \setminus E_0)\cap E(R)|$ is even by Lemma~\ref{R-lem}~(b).
 
 Let $H$ be the subgraph induced by the edge set $E_G(C') \cap E_G(X', Y')\cap E(R)=(E_G(C') \setminus E_0)\cap E(R)$.
Since $X' \subseteq X\setminus V(\mathcal{P}_0)$, we can apply (\ref{X1X2}) to this edge set. Hence
\begin{itemize}
    \item for each vertex $y\in Y' \setminus S$, it is incident with exactly two edges in $H$. Additionally, due to the properties of $V_1$ and $V_2$, these two edges are either both in $E(R)$ or both in $E(G) \setminus E(R)$;
    \item for $u_i^1$ and $u_i^2$ with $i \in \{r+1, r+2, \dots, |Z'|\}$, both $u_i^1$ and $u_i^2$ have degree $1$ in $H$ and  are either both in $V_1$ or both in $V_2$. 
    Consequently, the two edges incident to them are  either both in $E(R)$ or both in $E(G) \setminus E(R)$;
    \item for $w_i^1$ and $w_i^2$ with $i \in \{1,2,\dots, |I|+|\mathcal{P}|\}$, both $w_i^1$ and $w_i^2$ have degree $1$ in $H$ and are contained in $V_2$. Thus the two edges incident to them are both in $E(R)$.
\end{itemize}
Therefore, $|E(H)|=|(E_G(C') \setminus E_0) \cap E(R)|$ is even.

Recall that $x$ and $y$ are chosen randomly. Hence, in the HC Construction of $C$, we can choose $y \in V_2$ such that $v_1^1y \in E(R)$ by (\ref{V2-Z}). This follows that $xy \in E(R)$ since $E_G(X\setminus E(\mathcal{P}_0),V_2)\subseteq E(R)$ by (\ref{X1X2}).
By above assertion, 
we conclude that $\big|E_G(\mathcal{P}) \cap E(R)\big|$ is even.

Finally, we are going to show
\[
    E_G\big(X \setminus \text{In}(\mathcal{P}), V_1\big) \subseteq E(G) \setminus E(R) \quad \text{and} \quad E_G\big(X \setminus \text{In}(\mathcal{P}), V_2\big) \subseteq E(R).
\]
Recall that (Claim~\ref{E(X,Y)} (ii))
\[
    E_G\big(X \setminus \text{In}(\mathcal{P}_0), V_1\big) \subseteq E(G) \setminus E(R) \quad \text{and} \quad E_G\big(X \setminus \text{In}(\mathcal{P}_0), V_2\big) \subseteq E(R),
\]
it is sufficient to show
\[
    E_G\big(\text{End}(\mathcal{P}) \cup I, V_1\big) \subseteq E(G) \setminus E(R) \quad \text{and} \quad E_G\big(\text{End}(\mathcal{P}) \cup I, V_2\big) \subseteq E(R),
\]
because $\text{In}(\mathcal{P}_0) \setminus \text{In}(\mathcal{P}) \subseteq V(\mathcal{P}_0) \setminus \text{In}(\mathcal{P}) \subseteq I \cup \text{End}(\mathcal{P})$.
% $X \setminus \text{In}(\mathcal{P}) \subseteq \big(X \setminus \text{In}(\mathcal{P}_0)\big) \cup \text{End}(\mathcal{P}) \cup I$.
Recall that $\text{End}(\mathcal{P}) \cup I=\bigcup_{i=1}^{|\mathcal{P}|+|I|}\{v_i^1, v_i^2\}$.

Arbitrarily choose an edge $e\in  E_G\big(\text{End}(\mathcal{P}) \cup I, Y\cup Z\big)$. Without loss of generality, assume $e=v_1^1z$.\footnote{Here it is reasonable to choose \( v_1^1 \) here. In HC Construction, we selected \( y \) from the neighbors of \( v_1^1 \). However, we can alternatively choose \( y \) from the neighborhood of the endpoint of \( e \) in \( X \).} In the HC construction we choose $y=z$.
Combining with our assertion and the fact $|E(\mathcal{P})\cap E(R)|$ is even, this follows that 
$|\{xz,zv_1^1\} \cap E(R)|$ is even.
Since $xz \in E(G)\setminus E(R)$ if $z \in V_1$ and $xz \in E(R)$ if $z \in V_2$, we conclude that
$v_1^1z \in E(G)\setminus E(R)$ if $z \in V_1$, and $v_1^1z \in E(R)$ if $z \in V_2$. This means that
$E_G\big(X\setminus \text{In}(\mathcal{P}), V_1\big) \subseteq E(G) \setminus E(R)$ and $ E_G\big(X\setminus \text{In}(\mathcal{P}), V_2\big) \subseteq E(R)$.

\end{proof}

           \begin{claim}\label{XVi}
               $E_G(X,V_1)\subseteq E(G)\setminus E(R)$, $E_G(X,V_2) \subseteq E(R)$, and $E_G(X)\subseteq E(G)\setminus E(R)$.
           \end{claim}
           \begin{proof}
           If \( t \leq -1 \), then \( |E(\mathcal{P}_0)| = 0 \). Therefore, by Claim~\ref{E(X,Y)}, we have 
             \[
                E_G(X, V_1) \subseteq E(G)\setminus E(R)  \quad \text{and} \quad E_G(X, V_2) \subseteq E(R).
             \]
           If \( E_G(X) \neq \emptyset \), arbitrarily choose an edge \( e \in E_G(X) \). Let \( \mathcal{P}=\{e\} \). Then $\mathcal{P}$ is a UDP in $G[X]$ with  \(E(\mathcal{P})=\{ e\} \). Hence by Claim~\ref{X-path}, we deduce that \( e \in E(G) \setminus E(R) \). 
           This implies that 
             \[
                E_G(X) \subseteq E(G) \setminus E(R).
             \]

           Now assume \( t \geq 1 \).  Then $E(\mathcal{P}_0)=t+1 \ge 2$. For any edge \( e = ab \in E(\mathcal{P}_0) \), define \( \mathcal{P}_e=\mathcal{P}_0 - e \). Then \( \mathcal{P}_e \) is a new UDP with \( t \) edges. 
         By  Claim~\ref{X-path}, we have  \( |E(\mathcal{P}_e) \cap E(R)| \) is even, and 
               \[
                  E_G\big(\{a,b\}, V_1\big) \subseteq E(G) \setminus E(R) \quad \text{and} \quad E_G\big(\{a,b\}, V_2\big) \subseteq E(R),
               \]
           since \( \{a, b\} \subseteq X\setminus \text{In}(\mathcal{P}_e) \). 
 %          Therefore, we can apply Claim~\ref{X-path} to \( \mathcal{P}(e) \). Consequently, \( |E(\mathcal{P}(e)) \cap E(R)| \) is even, and 
 %              \[
 %                  E_G\big(\{a,b\}, V_1\big) \subseteq E(G) \setminus E(R) \quad \text{and} \quad E_G\big(\{a,b\}, V_2\big) \subseteq E(R),
 %              \]
 %          since \( a, b \notin In(\mathcal{P}(e)) \). 
           Due to the random selection of \( e \) and the fact that \( |E(\mathcal{P}_e)\cap E(R)| = t \) is odd, it follows that
             \[
                E_G(X, V_1) \cap E(R) = \emptyset, \quad E_G(X, V_2) \subseteq E(R), \quad \text{and} \quad E(\mathcal{P}_0) \subseteq E(G) \setminus E(R).
             \]

             %Since \( \mathcal{P}_0 \) is chosen randomly, we conclude that \( E(\mathcal{P}_1) \subseteq E(G) \setminus E(R) \) also holds for any UDP \( \mathcal{P}_1 \) with \( \max\{t+1, 0\} \) edges.
             
We next demonstrate that $E_G(X)\subseteq E(G)\setminus E(R)$. For any \( e \in E_G(X)\),  define a UDP
$$\mathcal{P}_e = \mathcal{P}_0 - \{e_1, e_2\} + e, $$ where \( e_1, e_2 \in E(\mathcal{P}_0) \) satisfying that:
\begin{itemize}
    \item  if \( e \in E_G(X \setminus V(\mathcal{P}_0)) \), select \( e_1\) and \( e_2 \) arbitrarily in $\mathcal{P}_0$;
    \item if \( e \in E_G(V(\mathcal{P}_0), X \setminus V(\mathcal{P}_0)) \), choose one edge from \( \{e_1, e_2\} \) to be incident with \( e \);
    \item if \( e \in E_G(V(\mathcal{P}_0)) \setminus E_G(\mathcal{P}_0) \), select \( e_1\), \( e_2 \) such that each is incident with one endpoint of \( e \).
\end{itemize}
Then for any edge \( e \in E_G(X) \setminus E_G(\mathcal{P}_e) \), we have \( E_G(\mathcal{P}_e) \setminus \{e\} \subseteq E_G(\mathcal{P}_0) \subseteq E(G) \setminus E(R) \), and \( |E(\mathcal{P}_e)|=t \) is odd. By Claim~\ref{X-path}, \( |E(\mathcal{P}_e) \cap E(R)| \) is even. It follows that \( e \in E(G) \setminus E(R) \). Finally, combining the above consequently, we obtain that \( E_G(X) \subseteq E(G) \setminus E(R) \).
\end{proof}

From (\ref{X1X2}), we have 
\[
\big(E_G(V_1) \cup E_G(V_2)\big) \cap E(R) = \emptyset, \text{ and } E_G(V_1, V_2) \subseteq E(R).
\]
Furthermore, according to Claim~\ref{XVi}, we have 
\[
E(X, V_1) \cap E(R) = \emptyset, \quad E(X, V_2) \subseteq E(R), \quad E_G(X) \cap E(R) = \emptyset.
\]
We conclude that \( E(R) = E_G(X \cup V_1, V_2) \), a contradiction to the property (c) of $R$.

\end{proof}

\begin{lem}\label{loc:dense}
  Suppose $\frac{1}{n}\ll \gamma$.
  If an $n$-vertex graph $G$ is Hamilton-connected, has minimum degree $\delta(G) \ge \frac{n-1}{2}$, and satisfies $e(A, B) \ge \frac{\gamma}{2} n^2$ for all (not necessarily disjoint) subsets $A, B\subseteq V(G)$ with $|A|,|B|\ge \frac{1-\gamma}{2}n$,  then
     $\mathcal{C}_{n}(G)=\mathcal{C}(G)$. 
 \end{lem}
\begin{proof}
Suppose that $\mathcal{C}(G) \neq \mathcal{C}_n(G)$. According to Lemma \ref{R-lem}, there exists a subgraph $R$ of $G$ satisfying (a), (b), and (c).
From (c), we have $\delta(R)\ge\frac 12\delta(G)$.

\begin{claim} \label{den:con}
    $G$ is $(\frac{\gamma}{4} n,3)$-connected.
\end{claim}
\begin{proof}
    For any subset $U\subseteq V(G)$ of size at most $\frac{\gamma}4  n$ and $u,v \in V(G)\setminus U$, we have $|N_{G}(u)\setminus U|,|N_G(v)\setminus U| \ge \delta(G) - |{U}| \ge \frac{n-1}{2}-\frac{\gamma}{4} n \ge \frac{1-\gamma}{2}n$. Thus $e(N_{G}(u)\setminus U, N_G(v)\setminus U)\ge \frac{\gamma}{2} n^2$. Select an edge $u'v'\in E(N_{G}(u)\setminus U, N_G(v)\setminus U)$ with $u' \in N_{G}(u)\setminus U, v' \in N_G(v)\setminus U$. Then $P=uu'v'v$ is a $u$-$v$-path of length three in $G-U$. It follows that $G$ is $(\frac{\gamma}{4}n,3)$-connected.
\end{proof}

\begin{claim} \label{Ham:con}
  For any $K \subseteq V(G)$ with $|K| \ge (1-\frac{\gamma}4)n+3$, $G[K]$ is Hamilton-connected.
\end{claim}
\begin{proof}
Clearly, $|V(G)\setminus K|\le \frac{\gamma}4n-3$.  By Claim~\ref{den:con}, we see that $G[K]$ is connected. For any two vertices $u,v \in K$, choose a longest $u$-$v$-path  $P_{uv}$ in $G[K]$. 
According to the Dirac's Theorem~\cite{D52}, we have $|E(P_{uv})|\ge 2\delta(G[K])+1 \ge 2(\delta(G)-|V(G)\setminus K|)+1 \ge (1-\frac{\gamma}{2})n+6$.
Label $P_{uv}=v_1v_2\dots v_{\ell_0}$ with $v=v_1$ and $u=v_{\ell_0}$, and for any $u',v' \in V(P_{uv})$, denote the subpath of $P_{uv}$ connecting $u'$ and $v'$ as $P_{u'v'}$.
We claim that $V(P_{uv})=K$. Otherwise, for any $w \in K\setminus V(P_{uv})$, we know that $|N(w)\cap V(P_{uv})| \ge \delta(G) - |\overline{V(P_{uv})}|\ge \frac{1-\gamma}{2}n+5$. Let $N^+_{P_{uv}}(w)=\{v_i: 2 \le i \le \ell_0-1, v_{i+1} \in N(w)\cap V(P_{uv})\}$. Then $|N^+_{P_{uv}}(w)|\ge |N(w)\cap V(P_{uv})|-2 \ge \frac{1-\gamma}{2}n+3$.  By the assumption for $A=B=N^+_{P_{uv}}(w)$, we have
\begin{equation*} \label{outedge}
    e(A,B)=e(N^+_{P_{uv}}(w))\ge \frac{\gamma}{2}n^2 \ge 1+|E(P_{uv})|,
\end{equation*}
this implies that we can find an edge $e=v_{i_1-1}v_{i_2-1} \in E_G(N^+_{P_{uv}}(w))\setminus E(P_{uv})$ with $i_1 < i_2$.
%, which is nonempty by (\ref{outedge}). 
By the definition of $E_G(N^+_{P_{uv}}(w))$, we know $v_{i_1},v_{i_2} \in N(w)\cap V(P_{uv})$. Then we get a longer path %$(v,P_{vv_{i_1-1}},v_{i_1-1},v_{i_2-1}, P_{v_{i_2-1},v_{i_1}},v_{i_1},w,v_{i_2},P_{v_{i_2},u},u)$
$$P_{vv_{i_1-1}}+v_{i_1-1}v_{i_2-1}+P_{v_{i_2-1},v_{i_1}}+v_{i_1}wv_{i_2}+P_{v_{i_2}u}$$
between $u$ and $v$ in $G[K]$, a contradiction.
\end{proof}

\begin{claim}\label{CL: S}
    For any $S \subseteq V(G)$ of size $|S| \le 100$ and any two vertices $x, y \in V(G)\setminus S$, there is an $x$-$y$-path in $R-S$ of length at most $16$.
\end{claim}
\begin{proof}
    Let $G'=G-S$ and $R'=R-S$. Then we have $\delta(G') \ge \frac{n-1}{2}-100$ and $\delta(R') \ge\frac 12\delta(G)-100\ge\frac{n-1}{4}-100$.
   For any $v \in V(G')$, let $N^0(v)=\{v\}$ and $N^{i}(v)=N^{i-1}(v) \cup N_{R'}\big(N^{i-1}(v)\big)$ for $i\ge 1$. 
   Then 
   \begin{align}\label{N^i}
  |N^{i}(v)| \ge \delta(R')+1\ge \frac{n-1}{4}-99 \text{ for any $v \in V(G')$ and $i\ge 1$}.    
   \end{align}
%   As $N_{R'}(v) \subseteq N^{t}(v)$ {\color{red} we have $|N^{t}(v)| \ge \frac{n-1}{4}-100$ for any $v \in V(G')$.} 
    To prove that for $x, y \in V(G)\setminus S$, there is an $x$-$y$-path in $R-S$ of length at most $16$, it is sufficient to show that $N^{8}(x) \cap N^{8}(y) \ne \emptyset$. Assume that there are $x,y \in V(G)$ satisfying $N^{8}(x) \cap N^{8}(y) = \emptyset$. Then $N^{t}(x)\cap N^{t}(y)=\emptyset$ and $E_R(N^{t}(x), N^{t'}(y))=\emptyset$ for any $t, t'\in [7]$.
    
We first claim that there are no vertices $w,z \in V(G')\setminus \{x,y\}$ such that $N^{1}(w)$, $N^{1}(x)$, $N^{1}(y)$, $N^{1}(z)$ are pairwise disjoint and there is no edge of $E(R)$ between them. 
 Suppose it is not the case, denote $N^1(u,v)=N^1(u)\cup N^1(v)$ for $u,v\in\{w,x,y,z\}$. Then $|N^1(u,v)|=|N^1(u)|+|N^1(v)|\ge 2(\frac{n-1}4-99)=\frac{n-1}2-198$
    and  $e_R(N^{1}(w,x), N^{1}(y,z))=0$. Note that $|V(G)\setminus (N^{1}(w,x)\cup N^1(y,z))|\le 399$. 
    Hence \begin{align*}
        e_R(N^{1}(w,x),\overline{N^{1}(w,x)})&=e_R(N^{1}(w,x), V(G)\setminus (N^{1}(w,x)\cup N^1(y,z))) \le 400n < \frac{\gamma n^2}{2}.
    \end{align*}
   However, according to the assumption (c) and Lemma~\ref{Class} (1),
    $$
    e_R(N^{1}(w, x),\overline{N^{1}(w,x)}) \ge \frac 12e_G(N^{1}(w, x),\overline{N^{1}(w,x)}) \ge  
    \frac{\gamma n^2}{2},
    $$
    a contradiction. 
    
    Second, we claim that there is no vertex $z \in V(G')\setminus \{x,y\}$ such that $N^{3}(x)$, $N^{3}(y)$, $N^{3}(z)$ are pairwise disjoint and there is no edge of $E(R)$ between them. Suppose it is not the case. If there exists $w \in V(G')\setminus (N^{3}(x)\cup N^{3}(y)\cup N^{3}(z))$, then we deduce that $N^{1}(w)$, $N^{1}(x)$, $N^{1}(y)$, $N^{1}(z)$ are pairwise disjoint and there is no edge of $E(R)$ between them, this is a contradiction. Now assume that $N^{3}(x)\cup N^{3}(y)\cup N^{3}(z) = V(G'),$ and $| N^{3}(x) |$ is the smallest one among them. Then $|N^3(x)|\le \frac n3$.
    Thus again according to the assumption (c) and  (\ref{N^i}),
  \begin{align*}
  e_R(N^{3}(x),\overline{N^{3}(x)})&\ge \frac12e_G(N^{3}(x),\overline{N^{3}(x)}) \ge\frac12\left(\delta(G')-\frac n3\right)|N^3(x)|\\
                                   &\ge\frac12 \left(\frac{n}{6}-101\right)\left(\frac{n-1}4-99\right)>\frac{\gamma}2 n^2.
  \end{align*} 
   However, $e_R(N^{3}(x),\overline{N^{3}(x)})=e_{R}(N^3(x), V(G)\setminus V(G')) \le 100n < \frac{\gamma n^2}{2}$, leading to a contradiction.
    
   At last, we assert that it is impossible that $N^7(x)$ and $ N^7(y)$ are disjoint and there is no edge of $E(R)$ between them.
    If there exists $z \in V(G)\setminus (N^7(x)\cup N^7(y))$, then $N^{3}(x)$, $N^{3}(y)$, $N^{3}(z)$ are pairwise disjoint, this is a contradiction.
Therefore, $V(G')= N^7(x)\cup N^7(y)$. Since $G$ is $(\frac{\gamma}{4} n,3)$-connected and $|S| \le 100$, we have $G'$ is connected. Then again according to
    the assumption (c),
  \begin{align*}
  e_R(N^{7}(x),{N^{7}(y)})&\ge \frac12e_G(N^{7}(x),\overline{N^{7}(y)}) > 0.
  \end{align*} 
It concludes that $N^8(x)\cap N^8(y) \neq \emptyset$, leading to a contradiction.

\iffalse
    This implies that there exists a vertex $z \in Z$ such that $N^{3}(x) \cup N^{3}(y)\cup N^{3}(z) = V(G')$. Otherwise, if there is a vertex $w\in V(G')\setminus(N^{3}(x) \cup N^{3}(y)\cup N^{3}(z))$, then we have $w,z \in V(G')\setminus \{x,y\}$ such that $N^{1}(w)$, $N^{1}(x)$, $N^{1}(y)$, $N^{1}(z)$ are pairwise disjoint and there is no edge of $E(R)$ between them, leading to a contradiction.
    Suppose that $| N^{4}(x) |$ is the smallest one among them. Then $|N^4(x)|\le \frac n3$.
    If $ N^{4}(x), N^{4}(y),
    N^{4}(z)$ are pairwise disjoint and there is no edge of $E(R)$ between them, then again according to the assumption (c) and Lemma~\ref{Class} (1),
  \begin{align*}
  e_R(N^{4}(x),\overline{N^{4}(x)})&\ge \frac12e_G(N^{4}(x),\overline{N^{4}(x)}) \ge\frac12\left(\delta(G')-\frac n3\right)|N^4(x)|\\
                                   &\ge\frac12 \left(\frac{n}{6}-101\right)\left(\frac{n-1}4-99\right)>\frac{\gamma}2 n^2.
  \end{align*} 
   However, $e_R(N^{4}(x),\overline{N^{4}(x)})=e_{R}(N^4(x), V(G)\setminus V(G')) \le 100n < \frac{\gamma n^2}{2}$, leading to a contradiction.
   This implies that 
    $N^{13}(x) \cup N^{13}(y) = V(G')$. 
    
    Since $e_R(N^{13}(x),N^{13}(y))\ge \frac{1}{2}e_G(N^{13}(x),N^{13}(y))\ge \frac{\gamma n^2}{4}$, $N^{14}(x)\cap N^{14}(y) \ne \emptyset$. 
\fi
\end{proof}

According to Claim~\ref{CL: S}, we can apply Lemma \ref{Switcher} to $G$ with respect to $R$, 
then we have an even cycle $C=(v_1,v_2,...,v_{h}, v_{h+1},...,v_{2h}, v_1)$ containing an odd number of edges from $R$ with $h \le 17$.

 As $G$ is $(\frac{\gamma}{4} n,3)$-connected, then $G-\{v_1,v_{h+1}\}$ is $(\frac{\gamma}{4} n-2,3)$-connected. Note that $h-1 \le \left\lceil(\frac{\gamma n}{4}-2)/4\right\rceil$. Applying Lemma~\ref{bidiameter} to $G-\{v_1,v_{h+1}\}$ with respect to $h-1$ pairs of vertices $(v_2,v_{2h}),(v_3,v_{2h-1}),\dots,(v_{h},v_{h+2})$,  we obtain $h-1$ disjoint paths $P_2,P_3,\dots,P_h$ in $G-\{v_1,v_{h+1}\}$ such that $P_{i}$ connects $v_{i}$ and $v_{2h+2-i}$ and $e(P_i)\le 3$.
 Set $P_1=v_1$ and $P_{h+1}=v_{h+1}$.    
We construct an $R$-parity switcher $W=(R, C_{(v_1,v_2,\dots, v_{2h})}, P_1, P_2,\dots,P_h, P_{h+1}, v_1, v_{h+1})$ with $|V(W)|\le 4h$.
    Let $U=(V(G) \setminus V(W)) \cup \{v_1,v_{h+1}\}$. Then $|U|\ge n-4h+2\ge(1-\frac{\gamma}4)n$. According to Claim~\ref{Ham:con},  we can find a Hamilton path $P_1$ connecting $v_1$ and $v_{h+1}$ within $G[U]$, leading to a contradiction with Corollary~\ref{parity}.

\end{proof}

\section{Proof of Theorem \ref{Odd}}
\Odd*
\begin{proof}
We choose $\frac{1}{n} \ll \gamma \ll \frac{\gamma}{\eta} \ll \alpha \ll \beta \ll \eta \ll 1$ and $n$ is odd.
Let $G$ be a Hamilton-connected graph on $n$ vertices with minimum degree $\delta(G) \ge \frac{n-1}{2}$.
By applying Lemma~\ref{Class} to $G$ with the parameter $\frac{1}{n} \ll \gamma$, we conclude that at least one of the statements (1), (2), or (3) must hold.

\begin{case}\label{Case 1}
  If $e(A,B) \ge \frac{\gamma}{2}n^2$ for  all $A$ and $B$ with at least $\frac{1-\gamma}{2}n$ vertices (not necessarily disjoint) (Lemma~\ref{Class} (1)). 
\end{case}
 This case is established by Lemma~\ref{loc:dense}.
    
\begin{case}\label{case2}
    There exists a set $A\subseteq V(G)$ of size $(\frac 12-21\gamma)n \le |A| \le \left(\frac 12 + 21 \gamma \right)n$ such that $e(A,\overline{A}) \le 4\gamma n^2$, and $\min\{\delta(G[A]),\delta(G[\overline{A}])\}\ge \frac{n}{5}$ (Lemma~\ref{Class} (2)). 
\end{case}

%Firstly, we can pick out all vertices who have many neighbors in both $A$ and $\overline{A}$.
    %By symmetry of $A$ and $\overline{A}$, we assume that $\frac{n}{2} \le |A| \le \left(\frac 12 + 21 \gamma \right)n$.
    Let 
    \begin{equation*}\label{biclique-1}
        Z = \{x \in A:|N(x) \cap \overline{A}| \ge \eta n \} \cup \{x \in \overline{A}:|N(x) \cap {A}| \ge \eta n \},
    \end{equation*}
    $X=A\setminus Z$, and $Y=\overline{A}\setminus Z$.
   %{, and $M$ be a maximum matching in $G[X,Y]$}.  
    Then $\Delta_X(Y), \Delta_Y(X) \le \eta n$. 
 Furthermore,
    by
     \begin{equation*}\label{}
       4\gamma n^2 \ge e(A,\overline{A}) \ge \frac{1}{2}(e(Z \cap A, \overline{A})+e(Z \cap \overline{A},A)) \ge \frac{|Z| \eta n}2,
    \end{equation*}
we have the following bounds by simple calculations:
\begin{equation} \label{biclique-2}
    |Z| \le \frac{8 \gamma n}{\eta} \le \frac{\alpha n}{2},
    \text{ and } \left(\frac{1}{2}- \alpha\right)n \le |X|,|Y| \le \left(\frac{1}{2}+ \alpha\right)n.
\end{equation}
Therefore, $\delta(X)\ge \delta(G)-|Z|-\Delta_Y(X)\ge \left(\frac{1}{2}- 2\eta\right)n$, and similarly,  $\delta(Y)\ge \left(\frac{1}{2}- 2\eta\right)n$,  and $\min\{\delta_Y(Z), \delta_{X}(Z)\}\ge \eta n - |Z| \ge \frac{3}{4} \eta n$,  and
$ e(X,Y) \le e(A,\bar{A})\le \alpha n^2$.

We claim that inequality (\ref{cod:biclique}) in Lemma~\ref{biclique} holds. When $|Z| \ge 3$, the inequality holds trivially. For cases where $|Z| = i$ with $i \in [2]$,  the Hamilton-connectedness of $G$ guarantees  $m(X,Y)\ge 3-i$. Then we have $\frac{4}{3}|Z|+m(X,Y) \ge \frac{4}{3}i+3-i=\frac13 i +3 \ge \frac{10}{3}$. Now we assume that $|Z|=0$. Then $|X|+|Y|=n$. Hence
$$\delta_Y(X)+\delta_X(Y)\ge \delta(G)-(|X|-1)+\delta(G)-(|Y|-1)\ge (n-1)-(|X|+|Y|)+2=1.$$ 
Thus at least one of $\delta_Y(X)$ or $\delta_X(Y)$ must be at least one. Without loss of generality, assume $\delta_Y(X)\ge 1$. This implies $|E(G[X,Y])|\ge |X|=|A|\ge (\frac 12-21\gamma)n$. 
To bound the matching number, we invoke K\"{o}nig's theorem for bipartite graphs $H$, which states that the matching number $\nu(H)$ equals the vertex cover number $\beta(H)$.
Applying this to $G[X,Y]$, we have
$$\nu(G[X,Y])=\beta(G[X,Y])\ge \frac{|E(G[X,Y])|}{\Delta_Y(X)}\ge \frac{(\frac 12-21\gamma)n}{\eta n}\ge 4.$$
Therefore, the inequality (\ref{cod:biclique}) holds.

According to Lemma \ref{biclique}, $\mathcal{C}(G)=\mathcal{C}_{n}(G)$ holds in this case.

\begin{case}\label{case 3}
    There exists a set A of size $\left(\frac 12 - 25 \gamma\right)n \le |A| \le \left(\frac 12 + 25 \gamma\right)n$  such that 
    $e(A, \bar{A})\ge (\frac 14 - 5\gamma)n^2$,  $e(A) \le 6 \gamma n^2$, and $\delta(G[A,\bar{A}])\ge \frac n5$ (Lemma~\ref{Class} (3)). 
   % the bipartite graph induced by the edges between A and $\overline{A}$ has at least $(1/4 - 5\gamma)n^2$ edges and minimal degree $\frac n5$, and $e(A) \le 6 \gamma n^2$(Lemma~\ref{Class} (3)).
\end{case}
    Let 
    $$
        Z_1 =\left \{x \in A:|N(x) \cap A| \ge \eta n\right\},\,\,  Z_2=\left \{x \in \bar{A} : |N(x) \cap \bar{A}| \ge \eta n\right\},
    $$
   and
    $$
        Z_3=\left\{x \in \bar{A}:|N(x) \cap A| < \left(\frac12-\eta\right)n \right\}.
    $$
It is clear that \( Z_3 \subseteq Z_2 \) since \( \delta(G) \ge \frac{n-1}{2} \). Recall \( \gamma \ll \frac{\gamma}{\eta} \ll \alpha \ll \eta \). We claim that \( |Z_1|, |Z_3| \le \frac{1}{3} \alpha n \). 
First, from the inequality 
\[
e(Z_1) \le \binom{|Z_1|}{2} \le e(A) \le 6\gamma n^2,
\] 
we have \( |Z_1| \le \sqrt{12\gamma}n \le \frac{1}{3} \alpha n \). 
Next, combining the inequalities:
\[
|Z_3||A| - e_G(Z_3, A) \le |A||\overline{A}| - e_G(A, \overline{A}) \le \frac{n^2}{4} - \left(\frac{1}{4} - 5\gamma\right)n^2 = 5\gamma n^2,
\]
and
$$|Z_3||A| - e_G(Z_3, A) \ge \left(|A| - \left(\frac{1}{2} - \eta\right)n\right)|Z_3| \ge (\eta - 25\gamma)n|Z_3|,   $$
%\begin{align*}
% |Z_3||A| - e_G(Z_3, A) &\ge |Z_3||A| - \left(\frac{1}{2} - \eta\right)n|Z_3| \\
%& = \left(|A| - \left(\frac{1}{2} - \eta\right)n\right)|Z_3| \\
%&\ge (\eta - 25\gamma)n|Z_3|,   
%\end{align*}
we conclude that
\[
|Z_3| \le \frac{5\gamma}{\eta - 25\gamma}n \le \frac{1}{3} \alpha n.
\]

Define $Z=Z_1 \cup Z_2$ if $|Z_2| \le \frac{2}{3}\alpha n$, and define 
     $Z=Z_1 \cup Z_2' \cup Z_3$ if $|Z_2| > \frac{2}{3}\alpha n$, where $Z_2'$ is a subset of $Z_2 \setminus Z_3$ with $|Z_2'| = \frac{2}{3}\alpha n - |Z_3|$.
Then $|Z|\le \alpha n$.
 Let 
     $X=A\setminus Z$     and      $Y=\overline{A}\setminus Z$.
     Then $X\cup Y\cup Z$ forms a partition of $V(G)$. 
     
Based on the choice of $Z$, we have $\Delta(G[X]) \le \eta n$. Additionally, $\delta_Y(X) \ge \delta(G)-\eta n-|Z|\ge (\frac12-2\eta)n$, while
    $\delta_X(Y) \ge (\frac12-\eta)n-|Z|\ge (\frac12-2\eta)n$, as guaranteed by $Z_3 \cap Y = \emptyset$.
Furthermore,   both  $\delta_X(Z)$ and $\delta_Y(Z)$ are bounded below by $\eta n -|Z| \ge \frac{3}{4}\eta n$. The cardinality of $X$ satisfies $|X|=|A|-|Z_1| \ge (\frac12- \frac{\alpha}2)n$.
 If $|Z_2| \le \frac{2}3\alpha n$, then $Y= \overline{A} \setminus Z_2$. Hence we also have $\Delta(G[Y]) \le \eta n$ and $|Y|=|\bar{A}|-|Z_2| \ge \left(\frac12 - 25 \gamma -\frac{2}{3}\alpha\right)n\ge \left(\frac12- \alpha\right)n$. Therefore, $X$ and $Y$ are symmetry in this case.  Without loss of generality, we may assume that $|X| \ge |Y|$.
  If  $|Z_2| >\frac{2}3\alpha n$, then $Y =\overline{A}\setminus (Z_2'\cup Z_3)$. According to the definition of $Z_2'$, we have $|Z_2'\cup Z_3|=\frac{2}{3}\alpha n$. Thus
    $$
    \left(\frac{1}{2}- \alpha\right)n \le |Y| = |\overline{A}|-\frac{2}{3}\alpha n \le \left(\frac12 + 25 \gamma -\frac{2}{3}\alpha\right)n \le \left(\frac12- \frac{\alpha}2\right)n \le |X|.
    $$
In summary, the cardinalities of $X$, $Y$ and $Z$ are bounded below by 
\begin{equation}\label{bipartite-2}
(\frac{1}{2}- \alpha)n \le |Y| \le |X| \le (\frac{1}{2}+ \alpha)n, \text{ and } |Z| \le \alpha n.
\end{equation}
% Recall that $f(X)\le |X|$ and 
\begin{claim}\label{f(X)}
    $f(X) \ge 2\delta(X)$.
\end{claim}
\begin{proof}
    Recall $ f(X) = \max\{e(\mathcal{P}): \mathcal{P} \text{ is  a  UDP in } G[X]\}$. If $f(X) \ge 2\eta n$, then the claim holds trivially since $\delta(X)\le\Delta(X) \le \eta n$. Hence we assume that $f(X) \le 2\eta n$. Choose $\mathcal{P}$ as a maximal union of paths in $G[X]$ (in the sense of number of edges). 

    Since $|V(\mathcal{P})| \le f(X) \le 2\eta n$, we can choose a vertex $u \in X \setminus V(\mathcal{P})$. By the maximality of $|E(\mathcal{P})|$, we assert that $N_G(u) \cap X \subseteq V(\mathcal{P})$. Moreover, any endpoint or any two consecutive vertices (of some path) in $\mathcal{P}$ should not be incident with $u$ in $G$, because otherwise, we could find a larger UDP by including $u$. Therefore, for any path $P$ in $\mathcal{P}$, we have $|N_G(u) \cap V(P)| \le \frac{|E(P)|}{2}$. This leads to the following chain of inequalities:
$$\delta(X) \le |N_G(u)\cap X| \le \sum_{P \in \mathcal{P}}|N_G(u)\cap V(P)| \le \sum_{P \in \mathcal{P}}\frac{|E(P)|}{2} = \frac{E(\mathcal{P})}{2} \le \frac{f(X)}{2}.$$ As we desired.
    
    %{\color{red} and $u \in X\setminus V(\mathcal{P})$. Consider $N_G(u) \cap V(\mathcal{P})$, denoted as $\mathcal{P}(u)$. It is easy to see that any endpoint or any two consecutive vertices (of some path) in $\mathcal{P}$ should not be in $\mathcal{P}(u)$, otherwise we can find a larger union of paths by embedding $u$ in it.}
\end{proof}

%Since {\color{red}  by checking the neighbors of vertices of a union of paths with $f(X)$ edges in $G[X]$,{\color{blue}Otherwise, we can find another union with more edges.}}
Combining Claim~\ref{f(X)} and the fact that $|X|+|Y|+|Z|=n$, we have $$\frac{n-1}2\le\delta(G) \le \delta(X)+|Y|+|Z| \le \frac{f(X)+n-(|X|-|Y|-|Z|)}{2}.$$
It follows that $|X|-|Y|-|Z|\le f(X)+1$.
%    where $f(X) \ge 2\delta(X) $ by checking the neighbors of vertices of a union of paths with $f(X)$ edges in $G[X]$.

If $|X|-|Y|-|Z|\ge f(X)$,     choose  two vertices $y_1,y_2 \in Y$. 
 Since $G$ is Hamilton-connected, there exists a Hamilton path $P_1=P_{y_1y_2}$ connecting $y_1$ and $y_2$. Let $\mathcal{P}_0$ be a UDP in $G[X]$ induced by all edges of $E(X) \cap E(\mathcal{P}_0)$. Then  $|E(\mathcal{P}_0)| \le f(X)$. 
 Since  each vertex $x \in X\setminus V({P}_1)$ satisfies $|N_{{P}_1}(x) \cap (Y\cup Z)| = 2$, while those in End$(\mathcal{P}_0)$ and in In$(\mathcal{P}_0)$ have $|N_{{P}_1}(x) \cap (Y\cup Z)|$ equal to 1 and 0, respectively,  we obtain \begin{align*}
 e_{{P}_1}(X,Y)&= 2(|X|-|V(\mathcal{P}_0)|)+|\text{End}(\mathcal{P}_0)|\\
               &=2|X|-2(|V(\mathcal{P}_0)|-\frac 12|\text{End}(\mathcal{P}_0)|)\\
               &=2|X|-2|E(\mathcal{P}_0)|\\
               &\ge 2|X|-2f(X) \\
               &\ge 2(|Y|+|Z|).     
 \end{align*}   
 %$$
 %   e_{{P}_1}(X,Y)= 2(|X|-|V(\mathcal{P}_0)|)+|End(\mathcal{P}_0)|
 %   \ge 2|X|-2f(X) \ge 2(|Y|+|Z|),
 %   $$
%    where $\frac12|End(\mathcal{Q})| + |E(\mathcal{Q})|=|V(\mathcal{Q})|$ holds for any union of paths $\mathcal{Q}$. 
 On the other hand, as $|N_{\mathcal{P}_1}(y) \cap (X\cup Z)| \le 2$ for any $y \in (Y\cup Z)\setminus \{y_1,y_2\}$ and $|N_{\mathcal{P}_1}(y) \cap (X\cup Z)| = 1$ for any $y \in \{y_1,y_2\}$, we have
    $$
    e_{P_1}(X,Y) \le 2(|Y|+|Z|-2)+2=2(|Y|+|Z|-1),
    $$
    we get a contradiction.
       So we have $f(X)\ge |X|-|Y|-|Z|+1$.
    By Lemma \ref{bipartite} we get $\mathcal{C}(G)=\mathcal{C}_{n}(G)$.
   
\end{proof}

\section{Concluding remarks}
In this paper, we  prove that for sufficiently large odd $n$, every Hamilton-connected graph $G$ on $n$ vertices with minimum degree at least $(n-1)/2$ is Hamilton-generated. Furthermore, the minimum degree and Hamilton-connected condition is tight to insure a graph is Hamilton-generated, This result completely answers Problem 1.3 aﬀirmatively. 

Note that if $G$ contains an odd cycle, a necessary condition for $C_n(G) = C(G)$ is that $G$ has an odd number of vertices.
When \( |V(G)| = n \) is an even number, any odd cycle in \( G \) cannot be represented as the symmetric difference of Hamiltonian cycles. 
This naturally motivates the following question.
\begin{prob}
  Let $n$ be even and large enough and let $G$ be an $n$-vertex graph. Whether the Dirac-type minimum degree $\delta(G)>n/2$ could  enable that $\mathcal{C}(G)=\mathcal{C}_{n}(G)+\mathcal{C}_{n-1}(G)$?
\end{prob}

Another interesting problem posed by Heinig is the following.
\begin{conj}[Heinig, \cite{Heinig}]
    If $n$ is large enough and $G$ is a balanced bipartite graph with $2n$ vertices and $\delta(G) \ge \frac n2 +1$, then $G$ is Hamilton-generated
\end{conj}

%Furthermore, the Hamilton-generated property of random graphs is an intriguing topic that warrants further exploration. Heinig offered a lower bound of $p$:

%\begin{thm}[Heinig, \cite{Heinig2}]
%  If $\varepsilon > 0$, $p \in [0,1]^{\mathbb{N}}$ with $p(n) \geq n^{-\frac{1}{2}+\varepsilon}$, then for $n \to \infty$ a random graph $G \sim G_{n,p}$ asymptotically almost surely has the following properties:

%\begin{enumerate}
%    \item $\mathcal{C}(G)=\mathcal{C}_{n}(G)+\mathcal{C}_{n-1}(G)$,
%    \item if $n$ is odd, then $\mathcal{C}_{n}(G)=\mathcal{C}(G)$.
%\end{enumerate}
%\end{thm}

%And the authors of \cite{Nenadov} improve this bound as follows, which is tight asymptotically.
%\begin{thm}[\cite{Nenadov}]\label{random}
%There exists \( C > 1 \) such that for \( p \geq C \log n / n \), 
%\( G \sim G(n, p) \) w.h.p.\ satisfies \( \mathcal{C}_n(G) = \mathcal{C}(G) \) if \( n \) is odd.
%\end{thm}

%Furthermore, Heinig \cite{Heinig2} showed that if $\mathcal{C}_n(G_{n,p}) = \mathcal{C}(G_{n,p})$ and $G(n,p)$ is not a forest, then necessarily $\delta(G_{n,p}) \geq 3$. 
%In \cite{Nenadov}, the authors proposed the following \textit{hitting time} problem, which would further refine Theorem \ref{random}:
%\begin{prob}[\cite{Nenadov}]
%    Suppose $n$ is odd, and consider the random graph process $\{G_m\}_{m \in (\frac{n}{2})}$. Is it true that, with high probability, $\delta(G_m) \geq 3$ implies $C_n(G_m) = C(G_m)$?
%\end{prob}

\end{document}